\documentclass[11pt]{amsart}
\usepackage{amsfonts,epsfig,fancyhdr}
\usepackage{newlfont,amsfonts,amssymb,amsmath,amsthm,amsgen,amscd,datetime,dsfont,hyperref,tensor}
\usepackage{multicol,picinpar,enumerate}

\usepackage{euscript,color}

\DeclareMathOperator{\hess}{\mathrm{Hess}}

\renewcommand{\O}{\ensuremath{\mathrm{O}}}

\newcommand{\vect}[1]{\boldsymbol{#1}}

\newcommand{\bmat}{\begin{pmatrix}}
\newcommand{\emat}{\end{pmatrix}}

\newcommand{\e}{\mathbf{e}}
\renewcommand{\d}{{\mathrm d}}
\newcommand{\bcase}{\begin{case}}
\newcommand{\ecase}{\end{case}}

\newcommand{\bclaim}{\begin{claim}}
\newcommand{\eclaim}{\end{claim}}

\newcommand{\bstep}{\begin{step}}
\newcommand{\estep}{\end{step}}

\newcommand{\bhlem}{\begin{hlem}}
\newcommand{\ehlem}{\end{hlem}}

\newcommand{\bleer}{\begin{leer}}
\newcommand{\eleer}{\end{leer}}
\newcommand{\bde}{\begin{de}}
\newcommand{\ede}{\end{de}}

\newcommand{\bs}{\begin{satz}}
\newcommand{\es}{\end{satz}}
\newcommand{\btheo}{\begin{theo}}
\newcommand{\etheo}{\end{theo}}
\newcommand{\bfolg}{\begin{folg}}
\newcommand{\efolg}{\end{folg}}
\newcommand{\blem}{\begin{lem}}
\newcommand{\elem}{\end{lem}}
\newcommand{\bnote}{\begin{note}}
\newcommand{\enote}{\end{note}}
\newcommand{\bprf}{\begin{proof}}
\newcommand{\eprf}{\end{proof}}
\newcommand{\bd}{\begin{displaymath}}
\newcommand{\ed}{\end{displaymath}}
\newcommand{\be}{\begin{eqnarray*}}
\newcommand{\ee}{\end{eqnarray*}}
\newcommand{\eeqa}{\end{eqnarray}}
\newcommand{\beqa}{\begin{eqnarray}}
\newcommand{\bi}{\begin{itemize}}
\newcommand{\ei}{\end{itemize}}
\newcommand{\bnum}{\begin{enumerate}}
\newcommand{\enum}{\end{enumerate}}

\newcommand{\beq}{\begin{equation}}
\newcommand{\eeq}{\end{equation}}

\newcommand{\rr}{\mathds{R}}

\newcommand{\M}{{\cal M}}

\newcommand{\vf}{\varphi}
\newcommand{\earr}{\end{array}\]}
\newcommand{\barr}{\[\begin{array}}
\newcommand{\bvec}{\left(\begin{array}{c}}
\newcommand{\evec}{\end{array}\right)}

\newcommand{\sumi}{\sum_{i=1}^n}

\newcommand{\g}{\mathfrak{g}}
\renewcommand{\a}{\mathfrak{a}}
\newcommand{\he}{\mathfrak{he}}
\newcommand{\siml}{\mathfrak{sim}}

\newcommand{\h}{\mathfrak{h}}

\newcommand{\m}{\mathfrak{m}}

\renewcommand{\k}{\mathfrak{k}}

\newcommand{\stab}{\mathfrak{stab}}

\newcommand{\+}{\oplus}

\newcommand{\oo}{\vect{0}}
\newcommand{\x}{\vect{x}}

\newcommand{\rrn}{\mathds{R}^n}
\newcommand{\so}{\mathfrak{so}}
\newcommand{\co}{\mathfrak{co}}

\newcommand{\del}{\partial}

\newcommand{\bbem}{\begin{bem}}
\newcommand{\ebem}{\end{bem}}
\newcommand{\bbez}{\begin{bez}}
\newcommand{\ebez}{\end{bez}}
\newcommand{\bbsp}{\begin{bsp}}
\newcommand{\ebsp}{\end{bsp}}

\DeclareMathOperator{\grad}{\mathrm{grad}}
\newcommand{\wt}{\widetilde}

\newcommand{\nabk}{\nabla^{\cal K}}

\renewcommand{\gg}{\mathrm{g}}
\newcommand{\hh}{\mathrm{h}}

\newcommand{\RR}{{\mathrm{R}}}
\newcommand{\R}{{\mathrm{R}}}

\renewcommand{\epsilon}{\varepsilon}

\newcommand{\inter}{\makebox[7pt]{\rule{6pt}{.3pt}\rule{.3pt}{5pt}}\,}

\newcommand{\belabel}[1]{\begin{equation}\label{#1}}

\newtheorem{mthm}{Theorem}
\newtheorem{mcor}{Corollary}

\theoremstyle{definition}
\newtheorem{de}{Definition}[section]
\newtheorem{bem}[de]{Remark}
\newtheorem{bez}[de]{Notation}
\newtheorem{bsp}[de]{Example}
\newtheorem{conj}[de]{Conjecture}

\newtheorem*{bsp*}{Example}
\newtheorem*{def*}{Definition}
\theoremstyle{plain}
\newtheorem{lem}[de]{Lemma}
\newtheorem*{lem*}{Lemma}
\newtheorem{satz}[de]{Proposition}
\newtheorem{folg}[de]{Corollary}
\newtheorem{theo}[de]{Theorem}
\newtheorem*{theo*}{Theorem}
\newtheorem*{conj*}{Conjecture}

\numberwithin{equation}{section}

\begin{document}

\title
{Locally homogeneous pp-waves}

\thanks{This work was supported by the  Australian Research
Council via  the grants FT110100429 and DP120104582.}

%\date{\today\ at \xxivtime}

\author{Wolfgang Globke}
\author{Thomas Leistner}\address{School of Mathematical Sciences, University of Adelaide, SA 5005, Australia}\email{wolfgang.globke@adelaide.edu.au, thomas.leistner@adelaide.edu.au}

\begin{abstract}
We show that every $n$-dimensional  locally homogeneous pp-wave is a plane wave, provided it is indecomposable and its curvature
operator, when acting on $2$-forms, has rank greater than one.
As a consequence we obtain that indecomposable, Ricci-flat  locally homogeneous pp-waves are plane waves. This generalises a
classical result by Jordan, Ehlers and Kundt 
in dimension $4$. Several examples  show that our assumptions on indecomposability and the rank of the curvature are essential.
\end{abstract}

\maketitle

%\tableofcontents

\section{Background  and  main results}
\label{intro}
A semi-Riemannian manifold $(\M,\gg)$ is {\em homogeneous} if it admits a transitive action by a group of isometries. This means that for each pair of points $p$ and $q$ in $\M$ there is an isometry of $(\M,\gg)$ that maps $p$ to $q$. In the spirit of Felix Klein's  {\em Erlanger Programm} to characterise geometries by their symmetry group, homogeneous manifolds are fundamental building blocks in geometry.
Homogeneity is strongly tied to the geometry and the curvature of a manifold. For example, homogeneous Riemannian manifolds are geodescially complete, and, as an example for the link to  curvature, we recall the celebrated result that any Ricci-flat homogeneous Riemannian manifold is flat
 \cite{alekseevsky-kimelfeld75}. A weaker version of homogeneity which still  guarantees  that the manifold looks the same everywhere is   local homogeneity: a semi-Riemannian manifold is {\em locally homogeneous} if   for each pair of points $p$ and $q$ in $\M$ there is an isometry  defined on a neighbourhood of $p$ that maps $p$ to $q$.

Here we will study  local homogeneity for a certain class of Lorentzian manifolds, the so-called {\em pp-waves} and the {\em plane waves}. 
Locally, an $(n+2)$-dimensional  {\em pp-wave} admits coordinates $(x^-,x^1, \ldots , x^n, x^+)$ such that \belabel{ppcoordintro}
\gg:=2\d x^+ (\d x^-+H\d x^+) +\delta_{ij}\d x^i \d x^j,
\end{equation}
where 
 $H=H(x^1, \ldots , x^n,x^+)$ is a function  not depending on $x^-$. For a {\em plane wave}, this function is required to be quadratic in the $x^i$'s with $x^+$-dependent coefficients. In general, they are not homogeneous, but they admit a parallel null (i.e.~non-zero and light-like) vector field. An invariant definition of pp-waves and plane waves is given as follows: 
% \bde\label{ppdef}
A Lorentzian manifold $(\cal M,\gg)$ 
is 
 a {\em pp-wave} if it admits a
 parallel null vector field $V \in \Gamma(T\cal M)$, i.e., $V\neq0$, $\gg(V,V)=0$ and $\nabla V=0$, and if  its curvature endomorphism $\R: \Lambda^2T\M\to  \Lambda^2T\M$ is non-zero and  satisfies
\begin{equation}\label{screen-flat1}
	\R|_{V^\bot\wedge V^\bot}=0, 
 \end{equation}
 where $V^\perp:=\{X\in T\M\mid \gg(X,V)=0\}$. 
%Furthermore, we say that   $(\M,\gg)$ is a {\em pp-wave in standard form} or a {\em standard pp-wave} if $\M$ is an open subset  of $\rr^{n+2}$ with global coordinates $(x^-,x^1, \ldots , x^n, x^+)$ and $H=H(x^1, \ldots , x^n,x^+)$ a function on $\M$ not depending on $x^-$, and $\gg$ is the Lorentzian metric defined by $H$,
%\belabel{ppcoord}
%\gg:=2dx^+\circ (dx^-+Hdx^+) +\delta_{ij}dx^i\circ dx^j,
%\end{equation}
 A {\em plane wave}  is a pp-wave with the additional condition
\belabel{planewave}
\nabla_U\R=0\quad \text{ for all }U\in V^\perp.
\eeq
%\ede
Four-dimensional pp-waves were  discovered in a mathematical
context by Brinkmann \cite{brinkmann25} as one class of Einstein
spaces that can be mapped conformally onto each other.
In physics,   plane waves and pp-waves  appeared in general relativity  \cite{EinsteinRosen37}, where they continue to play an important role (see for example  \cite {BondiPiraniRobinson59,jek60} for more references) as metrics for which  the  Einstein  equations become linear and, when they solve these equations,  describe the propagation of gravitational waves with  flat surfaces as wave fronts.
Later Penrose discovered that when ``zooming in on null geodesics'' every space-times has a plane wave as limit \cite{penrose76}. More recently, the conditions under which the homogeneity of a Lorentzian manifold is inherited by its Penrose limit were studied 
extensively by Figueroa-O'Farrill,
Meessen and Philip~\cite{Figueroa-OFarriPhilipMeessen05,Philip06}.  Moreover, having linear Einstein equations and  a large number of parallel spinor fields, higher-dimensional plane waves and pp-waves recently appeared as supergravity backgrounds, e.g. in \cite{hull84pp}, and there is now a vast amount of literature on them. For more recent results on homogeneity see the work by Figueroa-O'Farrill et al.~ in \cite{Figueroa-OFarriMeessenPhilip05, Figueroa-OFarriHustler12,Figueroa-OFarriHustler14}.

A systematic study of $4$-dimensional pp-waves was undertaken by Jordan, Ehlers and Kundt in \cite{jek60} (see  the English republication \cite{jek60engl} and also \cite{ehlers-kundt62}, 
 where   the name {\em pp-wave} for {\em plane fronted with parallel rays} was introduced). Among other aspects, in \cite{jek60}   the isometries of $4$-dimensional, gravitational (i.e.~Ricci-flat) pp-waves are considered and the   Killing equation is solved completely. As a consequence, the possible dimensions of the space of Killing vector fields are given and in each case the form of the metric is determined explicitly. This rather satisfying result allows \cite{jek60} to conclude:
\bnum[(A)]
\item \label{A} If  a $4$-dimensional  Ricci-flat pp-wave $(\M^4,\gg)$ is locally $V^\perp$-homogeneous, then it  is a plane wave.
In particular, if $(\M^4,\gg)$  is Ricci-flat and  locally homogeneous,
then it is a  plane wave.
\enum
Here,  {\em local $V^\perp$-homogeneity} is a generalisation of local homogeneity taking into account the parallel null vector field $V$ that exists on a pp-wave:  the distribution $V^\perp$ is parallel as well and defines a foliation of $\M$ into totally geodesic leaves of codimension $1$. Then we say that $(\M,\gg)$ is {\em locally $V^\perp$-homogeneous} if  for all pairs $p,q\in\M$ that are  {\em in the same leaf} of $V^\perp$, there is a neighbourhood $\cal U$ of $p$ in $\M$ and an isometry $\phi$  between $(\cal U,\gg)$ and $(\phi(\cal U),\gg)$ that maps $p$ to $q$.

Note that proving \eqref{A}  amounts to showing that local homogeneity (in $V^\perp$-directions) forces all third derivatives $\del_i\del_j\del_k H$, for $i,j,k=1, \ldots ,n$ to vanish. This is a much harder problem in  dimensions higher than $4$. The 
 methods used in \cite{jek60} in order to solve the Killing equation  are restricted to dimension $4$ and also use that the function $H$ is harmonic, as a consequence of Ricci-flatness.

%Most remarkably, the implication \eqref{B} is no longer true if the assumption of Ricci-flatness is dropped. A counter example was provided in \cite[Table II, row 9]{sippel-goenner86} (see also our Remark \ref{sgremark}), where the Killing vector fields were determined without the assumption of Ricci-flatness.
Statement \eqref{A}  is no longer true without the assumption of Ricci-flatness: Sippel and Goenner in \cite{sippel-goenner86} solved   the Killing equation for a $4$-dimensional pp-wave $(\M^4,\gg)$ without assuming $\mathrm{Ric}=0$ and gave an example of a homogeneous pp-wave  that is not a plane wave (see our Example \ref{sgex}). However, it turns out that the metric in this example decomposes into a product of a $3$-dimensional pp-wave and $\rr$. Note that in \cite{jek60} such a decomposition was impli\-citly excluded by the Ricci-flatness: if a $4$-dimensional Ricci-flat  manifolds splits as a Riemannian product, then it is flat.  Hence, the results in \cite[Table II, p. 1234]{sippel-goenner86} establish the following result:
\bnum[(B)]
\item \label{B} If  a $4$-dimensional  indecomposable  pp-wave $(\M^4,\gg)$ is locally $V^\perp$-homogeneous, then it  is a plane wave.
In particular, if $(\M^4,\gg)$  is indecomposable and   locally homogeneous,
then it is a  plane wave.
\enum
Here, when saying  
that the manifold is {\em indecomposable}, we mean that 
the holonomy algebra acts indecomposably.
Therefore, when looking for a generalisation of \eqref{A} or (B) to arbitrary dimensions the notion of indecomposability is relevant. We say that a semi-Riemaniann manifold $(\M,\gg)$ is {\em strongly indecomposable} if $(\M,\gg)$ does not split
as a local semi-Riemannian product anywhere, i.e, there is no point in $\M$ that has a neighbourhood on which $\gg$ is a product metric. Clearly, by 
the local version of the de Rham-Wu splitting theorem,
 the holonomy algebra of  a strongly indecomposable manifold acts indecomposably (i.e.~without non-degenerate invariant subspace), but the converse in general is not true. In addition to strong indecomposability we will need another condition on the curvature tensor $\R$ of a pp-wave. From the very definition of a pp-wave  it follows that the rank of $\R$ when acting on $2$-forms does not exceed $\dim(\M)-2$. 
When proving a generalisation of statement~(B), we have to assume that generically the rank of $\R$ is larger than $1$:
\begin{mthm}\label{maintheo}
Let $(\M,\gg)$ be a  pp-wave of arbitrary dimension  with parallel null vector field $V$. Assume that $(\M,\gg)$ is strongly indecomposable and in addition that  almost everywhere the rank of its curvature endomorphism acting on $\Lambda^2T\M$ is larger than one.
Then $(\M,\gg)$ is a plane wave
if it 
is locally $V^\perp$-homogeneous. 
%
% if for each $p\in \M$ there are
%%exist a neighbourhood $\cal U$ with 
%Killing vector fields  spanning 
%$V^\perp|_p$ at $p$.
\end{mthm}
Here by ``almost everywhere'' we mean that there is no open set on which the rank of the curvature endomorphism is $\le 1$.
Note that  the assumption on the rank of the curvature prevents us from applying Theorem \ref{maintheo} to  $3$-dimensional 
pp-waves, for which the rank of $\R$ cannot be bigger than $1$.
Indeed, in Example \ref{ex_dim3} we exhibit a $3$-dimensional, locally homogeneous pp-wave that is {\em not} a plane wave, which shows that the assumption on the rank is crucial. However, since  Ricci-flat pp-waves always satisfy this assumption  (see Lemma \ref{ricrank}), we obtain a generalisation of statement \eqref{A} to arbitrary dimensions:
%\begingroup
%\renewcommand\thefolg{\Alph{folg}}
% 
\begin{mcor}\label{folg1}
A strongly  indecomposable, Ricci-flat and locally $V^\perp$-homogeneous  pp-wave is a plane wave.
\end{mcor}

In locally homogeneous manifolds all points have isometric neighbourhoods.
Hence,  a locally homogeneous manifold is strongly indecomposable whenever it is indecomposable, and the rank of the curvature endomorphism is constant. This  yields

\begin{mcor}\label{folg2}
An indecomposable, locally homogeneous  pp-wave   is a plane wave if, at one point,  the rank of the curvature endomorphism is greater than one.
\end{mcor}

\begin{mcor}\label{folg3}
Indecomposable, Ricci-flat and locally homogeneous  pp-waves are  plane waves.
\end{mcor}
Corollary \ref{folg3} is an instance of the phenomenon that Ricci-flat pp-waves with some additional geometric conditions  have to be plane waves. Another instance of this phenomenon is given in \cite{leistner-schliebner13}, where it is shown that compact Ricci-flat pp-waves are plane waves.

When proving Theorem \ref{maintheo} we use the following property which is implied by local $V^\perp$-homogeneity (see Lemma~\ref{homkill}): for every $p\in\M$ there are Killing vector fields, defined on a neighbourhood of $p$, which, when evaluated at $p$, span $V^\perp|_p$. However, these Killing vector fields might not be local sections of $V^\perp$. If they are, we can drop the assumption on the rank of the curvature and obtain

\begin{mthm}\label{kVtheo}
Let $(\M,\gg)$ be a strongly indecomposable pp-wave in which each point admits a neighbourhood $\cal U$ with local Killing vector fields that span $V^\perp|_{\cal U}$. Then $(\M,\gg)$ is a plane wave.
\end{mthm}
This is a version of a result for {\em commuting} Killing vector fields tangent to $V^\perp$:
\begin{mthm}\label{jek}
Let $(\M,\gg)$ be a semi-Riemannian manifold of dimension $m$ and 
assume that  there are {\em commuting}  Killing vector fields that
 span a null distribution (i.e., a distribution on which $\gg$ degenerates) of rank $m-1$. Then
 $(\M,\gg)$ admits a parallel null vector field $V$ and its curvature satisfies 
 \[\R(X,Y)Z=0\ \text{ and }\  \nabla_X\R=0,\]
 for all $X,Y,Z\in V^\perp$. In particular, if $(\M,\gg)$ is Lorentzian, then it is a plane wave.
\end{mthm}
Jordan, Ehlers and Kundt
\cite[Theorem 4.5.2]{jek60} proved Theorem \ref{jek} for $4$-dimensional 
 Lorentzian manifolds, but their proof works in any dimension and signature (see our Section \ref{wavesection}). In contrast, our proofs of Theorems \ref{maintheo} and \ref{kVtheo} use completely different methods than those  in \cite{jek60}.
In fact,
our proof of Theorem \ref{maintheo} does not require a full solution of the Killing equation (which we derive in Section \ref{killsec})  but  a detailed analysis of its consequences (in Section~\ref{proofsec}).  Moreover, we use algebraic results such as the classification of subalgebras of the Lie algebra of  similarity transformations $\siml(n)=(\rr\+\so(n))\ltimes \rrn$ 
of $\rrn$
 that act indecomposably on $\rr^{1,n+1}$ via $\siml(n)\subset \so(1,n+1)$. This classification is due to B\'{e}rard-Bergery and Ikemakhen  \cite{bb-ike93},
%, which is due to Alekseevsky \cite{alekseevsky75, alekseevskij-vinberg-solodovnikov93}. 
%Using this result, our proof reveals an interesting relation between the Lie algebra of Killing vector fields and subalgebras of $\so(1,n+1)$ acting indecomposably but not irreducibly on $\rr^{1,n+1}$, which on the other 
and plays an important role in the classification of indecomposable Lorentzian holonomy algebras in \cite{leistnerjdg}.

As we have pointed out above,  Example \ref{ex_dim3} shows that, at least in  dimension $3$, the condition on the rank of the curvature cannot be dropped. However, obvious generalisations of Example \ref{ex_dim3} lead either to {\em non-homogeneous} pp-waves (as in \cite{ehlers-kundt62}, see our Example \ref{transex}) or to {\em decomposable} homogeneous pp-waves (as  in  \cite{sippel-goenner86}, our Example~\ref{sgex}). Hence, with statement (B) in mind, we are tempted to conjecture (see also Section \ref{transverseremark} for more details):
\begin{conj*} Any indecomposable locally homogeneous pp-wave of dimension larger than $3$ is a plane wave.
\end{conj*}
In relation to this we should point out that the rank assumption is independent from the assumption of strong indecomposability: in Example \ref{rankexample} we present a $4$-dimensional, strongly indecomposable plane wave metric whose curvature has rank $1$. 

%\endgroup
Locally homogeneous plane waves turn out to be reductive (see Section \ref{pwred-sec}) and have been
been classified by Blau and O'Loughlin 
\cite{blau-oloughlin03} (see our Section \ref{rem_hom_PW}). As a consequence, 
with the exception of the curvature rank one case, our reduction to the plane waves   yields a classification 
 of indecomposable locally homogeneous pp-waves. The curvature rank one case remains open for further study. Also we believe that our methods employed in Section \ref{proofsec} are useful in a wider context and will give a better understanding of the more general class of indecomposable locally homogeneous Lorentzian manifolds.
  
 The paper is structured as follows: In Sections \ref{kill-sec}  we review facts about locally homogeneous spaces.   In Section \ref{wavesection} we present some facts about  pp-waves and plane waves including a fundamental coordinate description. In Section \ref{killsec} we  derive the Killing equation for pp-waves in these coordinates (Theorem \ref{killtheo}) and, moreover, use this to obtain some useful facts, including the reductivity of  homogeneous plane waves. We also review the classification of homogeneous plane waves in \cite{blau-oloughlin03} and  give a couple of examples that illustrate important features. Finally, in Section~\ref{proofsec} we will use the obtained results to prove Theorems \ref{maintheo} and \ref{kVtheo}. The appendix contains a proof of the coordinate description that turns out to be fundamental for our approach.
\subsection*{Acknowledgements} We would like to thank  Helga Baum for helpful discussions and Daniel Schliebner for commenting on the first version of this paper

\section{Killing vector fields and locally homogeneous spaces}
\label{kill-sec}
Let $(\M,\gg)$ be a semi-Riemannian manifold with Levi-Civita connection $\nabla$. A {\em Killing vector field} $K\in \Gamma (T\M)$ 
is a vector field whose flow $\phi_t$ consists of local isometries of $\gg$, i.e.~$\phi_t:(\cal U,\gg)\to (\phi_t(\cal U),\gg)$ is an isometry, where $\cal U$ is a neighbourhood of $p$ on which $\phi_t$ is defined. If $K$ is complete, then all $\phi_t$'s are global isometries.

Clearly, $K$ is a Killing vector field if and only if 
 the $(2,0)$-tensor 
$\gg(\nabla K,\cdot )$ is skew-symmetric, i.e.
\belabel{killingdef}
\gg(\nabla_X K,Y)+\gg(X,\nabla_Y K)=0
\quad\text{ for all }X,Y\in T\M.
\eeq
Let us denote the real vector space of Killing vector fields of $(\M,\gg)$ by $\k$. The Lie bracket of vector fields equips $\k$ with a Lie algebra structure.

In order to derive the integrability conditions for the Killing equation \eqref{killingdef}, we recall the classical approach by Kostant \cite{Kostant55}.
Let us denote by $\so(T\M,\gg):=\{\phi \in\mathrm{ End}(T\M)\mid \gg(\phi(X),Y)+\gg(\phi(Y),X) =0\}$ the bundle of skew-symmetric endomorphisms. For a 
 Killing vector field $K$, we  define the section $\phi^K:= \nabla K$ of $\so (T\M,\gg)$. A straightforward computation shows that the Killing equation \eqref{killingdef} implies that
\[
\nabla_X\phi^K = -\R(K,X),
\]
where $\R$ denotes the curvature tensor of $\gg$ defined by $\R(X,Y)=\left[\nabla_X,\nabla_Y\right]-\nabla_{[X,Y]}$. Hence, we  define the vector bundle 
\[\cal K:=T\M\+\so(T\M,\gg) \longrightarrow \M
\]
and furnish it with the covariant derivative
\[
\nabk_X
\begin{pmatrix}
K\\\phi
\end{pmatrix}
:=
\begin{pmatrix}
\nabla_X K -\phi(X)
\\
\nabla_X\phi +\R(K,X)
\end{pmatrix}.
\]
We get the vector space isomorphism
\[
\k\simeq \{\text{parallel sections of $(\cal K,\nabk)$}\},
\]
which shows that $\dim (\k)\leq \mathrm{rk}(\cal K)=\frac{1}{2}m(m+1)$, where $m=\dim(\M)$. It also shows that a Killing vector field $K$ is uniquely determined by  the values $K|_{p}\in T_p\M$ and $\nabla K|_{p}\in \so(T_p\M,\gg_p)$ at a point $p\in M$ and thus yields an injection of $\k$ into the Lie algebra of semi-Euclidean motions,
\belabel{eval}
\begin{array}{rcl}
\kappa\ :\ \k&\hookrightarrow &\so(r,s) \ltimes \rr^{r,s}\\
K&\mapsto & - \left (\left(\epsilon_i\epsilon_j\gg_p(\nabla_{\e_i}K,\e_j\right)_{i,j=1}^m,\left(\epsilon_k\gg_p( K_p,\e_k)\right)_{k=1}^m\right),
\end{array}
\eeq
where $(r,s)$ is the signature of $\gg$, $p\in \M$ and $\e_i$ an orthonormal basis of $T_p\M$, i.e., $\gg(\e_i,\e_j)=\epsilon_i\delta_{ij}$. 
Note the minus in front of the image. It ensures that for the flat metric on $\rr^{r,s}$ this map is a Lie algebra isomorphism (instead of an {\em anti}-isomorphism) between the Killing vector fields and the group of motions. In general,  this map is {\em not} a Lie algebra homomorphism. For example, the Killing vector fields of the $m$-sphere are isomorphic to $\so(m+1)$ rather than $\so(m)\ltimes \rr^n$.  
In fact, a lengthy but straightforward computation reveals
\[
\nabla [K,\hat K] =[\nabla K,\nabla\hat K] - \R(K,\hat K),
\]
where the right-hand side bracket is the commutator of linear maps, which yields
\belabel{lahom}
\kappa([K,\hat K])-\left[\kappa(K),\kappa(\hat K)\right]_{\so(t,s)\ltimes \rr^{t,s}} = -\left(\epsilon_i\epsilon_j \R_p(K_p,\hat K_p,\e_i,\e_j),0\right).
\eeq

%However, fixing a $p\in \M$ and  restricting $\kappa$ to the subalgebra
%\belabel{stabalg}\h(p):=\{K\in \k\mid K_p=0\},\eeq
%equation \eqref{lahom} shows that $\kappa$, when restricted to $\h(p)$,  becomes an injective Lie algebra homomorphism
%\marg{Don
%\[\kappa:\h(p)\hookrightarrow \so(t,s).\]

Returning to the integrability condition for the Killing equation, 
 we compute the curvature $\R^\cal K$ of $\nabk$ and we get
\[
\R^\cal K(X,Y)
\begin{pmatrix}
K\\\phi
\end{pmatrix}=
\begin{pmatrix}
0
\\(\nabla_K\R)(X,Y) - (\phi\cdot \R)(X,Y)
\end{pmatrix},
\]
where $ \phi\cdot \R$ denotes the canonical  action of an endomorphism on $(3,1)$-tensors. Hence
 the existence of a parallel section $(K,\phi)$ of $\cal K$ gives the integrability condition 
\belabel{kill-int}
\nabla_K\R=\phi\cdot \R,
\eeq
where $\phi=\nabla K$ and $\R$ is the curvature of $\gg$, i.e., we have for $X,Y,Z\in T\M$ that
\[
(\nabla_K\R)(X,Y)Z 
=
\phi (\R(X,Y)Z)-\R(\phi(X),Y)Z- \R(X,\phi(Y))Z- \R(X,Y)\phi(Z).
\]

Now, assume that $(\M,\gg)$ enjoys
the existence of a parallel vector field $V$. We define two vector spaces
\be
\k(V)&:=&\{K\in \k\mid \gg(K,V)=0\},\\
\k'(V)&:=& \{K\in \k\mid \nabla_VK=0\}
\ee
and observe
\blem\label{kill-lem}
If $V$ is a parallel vector field,  then we have the following inclusion of subalgebras
\[ \k(V)\subset \k'(V) \subset \k.\]
\elem
\bprf
First we check the inclusion $\k(V)\subset \k'(V)$. Indeed, for a Killing vector field $K\in \k$, the derivative of the function $\gg(V,K)$ satisfies
\belabel{evalV}X(\gg(K,V))=\gg(\nabla_XK,V)=-\gg(\nabla_VK,X).\eeq
First, this implies that if $K\in \k(V)$ then we also have $\nabla_VK=0$, i.e., $K\in \k'(V)$.

Next we note that both $\k(V)$ and $\k'(V)$ are subalgebras:
Clearly, if $V$ is parallel, $V^\perp$ is involutive and hence $\k(V)$ is closed under the bracket.
Moreover, for $K,\hat K\in \k'(V)$ we have that
\[
\nabla_V[K,\hat K]
=
\nabla_V\nabla_K\hat K-\nabla_V\nabla_{\hat K} K
=
\nabla_{[K,V]}\hat K-\nabla_{[\hat K,V]}K
=0,\]
since $V\inter \R=0$ and $[K,V]=-\nabla_VK=0$. Hence,  also $\k'(V)$ is a subalgebra.
\eprf
\blem
A Killing field $K$ satisfies $K\in\k(V)$ if and only if at some, and hence any, point $p\in \M$ we have
\belabel{killV} 
\gg(K,V)|_p=0, \quad \gg(\nabla_X K,V)|_p=0\ \text{ for all }X\in T_p\M.
\eeq
\elem

\bprf
First note that 
\eqref{evalV} implies that any $K\in \k(V)$ satisfies $\gg(\nabla_XK,V)\equiv 0$.

Conversely, assume $\gg(K,V)|_p=0$ and $\gg(\nabla_X K,V)|_p=0$ for all $X\in T_p\M$ at  $p\in \M$.
Let $\gamma $ be a geodesic emanating from $p$ with $\dot\gamma(0)=X$. Then by \eqref{evalV} we have
\be
\frac{\d^2}{\d t^2}(\gg(K,V)|_{\gamma(t)})
&=&
-\gg(\nabla_{\dot\gamma(t)}\nabla_V K,\dot\gamma(t)))
\\
&=&
-\gg(\nabla_V \nabla_{\dot\gamma(t)} K,\dot\gamma(t)))
- \gg(\nabla_{[\dot\gamma(t),V]} K,\dot\gamma(t)))
\\
&=&
\gg( \nabla_{\dot\gamma(t)} K,\nabla_V\dot\gamma(t)))
+ \gg(\nabla_{\dot\gamma(t)} K,[\dot\gamma(t),V]))
\\
&=&0,
\ee
i.e., $\gg(K,V)|_{\gamma(t)}$ is linear in $t$. Hence it is determined by its value and its derivative at $p$ which we both have assumed to be zero, forcing $\gg(K,V)|_{\gamma(t)}\equiv 0$. This shows that $\gg(K,V)$ is zero on a normal neighbourhood and thus zero everywhere.
\eprf
This lemma implies the following: Let $v\in\rr^{r,s}$  be the image of $V$ in $\rr^{r,s}$ under $\kappa$, i.e., $\kappa(V)=(0,v)$, and let $\stab(v)$ be its stabiliser in $\so(r,s)$. Then 
\[
\kappa:\k(V)\hookrightarrow \stab(v)\ltimes \rr^{r,s}.
\]
We will work with a different vector space of Killing vector fields, namely with
 \[
\k_p(V):=\{K\in \k\mid \gg(K,V)|_p=0\},
 \]
 for a fixed point $p\in \M$.
In general, this is not a Lie algebra. However, we will see that for pp-waves it is a Lie algebra, a fact which turns out to be very useful.

\bigskip

Now we consider locally homogeneous and locally $V^\perp$-homogeneous manifolds as defined in Section \ref{intro}. Both can be described in terms of Killing vectors.

 $(\M,\gg)$ is locally homogeneous if and only if for each point there exists Killing vector fields  spanning $T_p\M$ when evaluated at $p$, i.e., for each point, 
the evaluation map combined with the projection on $\rr^{r,s}$
\[\kappa:\k\to \so(r,s)\ltimes \rr^{r,s} \to \rr^{r,s}\]
is surjective. Moreover, $(\M,\gg)$ is homogeneous (the isometry group acts transitively on $\M$) if and only  if this holds for {\em complete} Killing vector fields. 

Analogously, we have
\blem\label{homkill}
Let $V$ be a parallel vector field on $(\M,\gg)$. If  $(\M,\gg)$  is locally $V^\perp$-homogeneous, then for each $p\in \M$ there exist local  Killing vector fields on a neighbourhood $\cal U$ of $p$ that span $V^\perp|_p$ when evaluated at $p$.
\elem
\bprf
Let $p\in \M$ and $X\in V^\perp|_p\subset T_p\M$. Let $\cal N_p$ be a leaf of $V^\perp$ through $p$  and $\xi:(-\epsilon, \epsilon) \to \cal N_p$ a curve such that $\dot\xi(0)=X$. Then, by assumption,  there is a curve $\gamma:(-\epsilon,\epsilon)\to G$ in the Lie pseudo-group  $G$ of local isometries around $p$ such that $\gamma_t(p)=\xi(t)$. Let $Y\in\g= T_{\mathrm{Id}}G$ be the tangent vector of $\gamma$ at $t=0$, i.e., $\dot \gamma(0)=Y$. Now let $\phi^Y_t:=\exp(tY)$ be the one-parameter pseudo-group that is defined by $Y$. This allows us to define a  vector field on a neighbourhood $\cal U$ of $p$ by 
\[
K(q):=\frac{\d}{\d t}( \phi_t^Y(q))|_{t=0},
\]
with $q\in \cal U$. Since the flow of $K$ is given by isometries, it is a Killing vector field. But we also have that $K(p)=X$ because
\[
K(p)=
\frac{\d}{\d t}( \exp(tY) (p))|_{t=0}
=
\d\Psi|_{\mathrm{Id}}\circ \d\exp|_0(Y)
=
\d\Psi|_{\mathrm{Id}}(Y)
=
\d\Psi|_{\mathrm{Id}}(\dot\gamma(0))
=
\dot \xi(0)=X,
\]
where $\Psi:G\to \M$ is defined by $\Psi(g)=g(p)$ and we use that $\d\exp|_0=\mathrm{Id}_{\g}$.
\eprf
 In general, these Killing vectors do not have to be tangent to $V^\perp$ everywhere.

Finally, note that a locally homogeneous manifold is strictly indecomposable (as defined in Section \ref{intro}) whenever it is indecomposable (i.e., the holonomy algebra acts indecomposably, that is without non-degenerate invariant subspace): if a locally homogenous manifold is a local product somewhere, it is a local product everywhere  and hence the holonomy algebra has a non-degenerate invariant subspace.

Unfortunately this does not hold in the case of local $V^\perp$-homogeneity for a parallel null vector field $V$. This can be easily seen for pp-waves  as in \eqref{ppcoordintro} on $\rr^{n+2}$: here the leaves of $V^\perp$ are given as $x^+=c$ constant. If $H(x^1, \ldots, x^n, x^+)\equiv 0$ for $x^+\in (a,b)$  but $\det(\del_i\del_j(H)|_{(x^1, \ldots, x^n,x^+)})\neq 0$ for some other  $x^+$, then the holonomy algebra acts indecomposably. However, near a point  with $x^+\in (a,b)$ the metric is flat.

\section{pp-waves and plane waves}
\label{wavesection}
Here we recall some basic properties of pp-waves and plane waves as defined in Section \ref{intro}.
First note that the defining equation  \eqref{screen-flat1} is equivalent to
\begin{equation}\label{screen-flat}
	%R|_{V^\bot\wedge V^\bot}=0,\ \text{ i.e., }\ 
	\R(U,W)=0\quad \text{ for all } U,W \in V^\bot,
\end{equation}
or to
\begin{equation}
	\label{screen-flat2}\R(X,Y)U \in \rr V\quad \text{ for all } U \in V^\bot \text{ and } X,Y\in T \cal M.
\end{equation}
A general pp-wave has an Abelian holonomy algebra contained in $\rrn$, where $\rrn$ is an Abelian ideal in the stabiliser $\so(n)\ltimes \rrn$  in $\so(1,n+1)$ of a null vector. The holonomy algebra is indecomposable if and only if it is equal to $\rrn$.
A pp-wave has the following coordinate description:

\blem\label{coordlemma}
Let $(\M,\gg)$ be a pp-wave and let $p\in \M$. Then there a are local coordinates $\vf=
(x^-,\x=(x^1, \ldots , x^n), x^+)$ on a neighbourhood $\cal U$ of $p$
and  a function $H\in C^\infty(\vf(\cal U))$ such that $H=H(x^+,\x)$ not depending on $x^-$
 such that, 
\belabel{ppcoord}
\gg=2\d x^+ (\d x^-+(H\circ \vf )\d x^+) +\delta_{ij}\d x^i \d x^j,
\end{equation}
where  $\delta_{ij}$ is the Kronecker symbol and where we use the summation convention.
In these coordinates the parallel null vector field is given by $V|_{\cal U}=\del_-:=\frac{\del}{\del x^-}$.
These coordinates are usually called {\em Brinkmann coordinates} after \cite{brinkmann25}. 

 Moreover, these coordinates can be chosen such that $\vf(p)=0$ and 
\belabel{delih0}
H(x^+, \oo)\equiv 0,\quad
\frac{\del H}{\del x^i}(x^+,\oo)\equiv 0,\end{equation}
for all $x^+$ from an interval around zero.
We call these coordinates {\em normal Brinkmann coordinates centred at $p$}.
\elem
Since the existence of coordinates as in \eqref{ppcoord} is well known, we only have to prove  normality, i.e., the property \eqref{delih0}. For sake of completeness we give a full proof of Lemma \ref{coordlemma} but 
defer it  to the appendix.

In Brinkmann coordinates
the non-vanishing components of $\nabla$ are
\belabel{nablai}
\begin{array}{rcl}
\nabla\del_i&=& \del_i(H)\d x^+\otimes \del_-
\\
\nabla\del_+&=& \d H\otimes \del_- -\d x^+\otimes  \grad(H),\end{array}\eeq
where $\grad(H)=\delta^{ij}\del_i(H)\del_j$ denotes the gradient of $H$ with respect to the flat metric $\delta_{ij}\d x^i \d x^j$ on $\rrn$.
This property justifies the term ``normal'' in Lemma \ref{coordlemma}: the covariant derivatives vanish at $\x=\oo$. The covariant derivatives of the corresponding
$1$-forms 
$\d x^i=\gg(\del_i,\cdot )$, $\d x^+=\gg(\del_-,\cdot )$ and $\d x^-=\gg(\del_+-2H\del_-,\cdot )$ are
\be
\nabla \d x^+&=&0
\\
\nabla \d x^i&=&\del_iH \d x^+\otimes \d x^+
\\
\nabla \d x^-&=&-2 \d H \d x^+.
\ee
For a pp-wave the parallel null vector field $V$ defines a parallel null distribution $V^\perp$ of rank $n+1$ for which the connection induced by the Levi-Civita connection  on the leaves of $V^\perp$  is flat. In Brinkmann coordinates, each leaf is defined by $x^+=c$ constant and  parametrised by the coordinates $x^-,x^1, \ldots , x^n$, and the formulae \eqref{nablai} show the flatness of the induced connection.
Moreover, equations \eqref{nablai} imply that  all the curvature components vanish apart from
\belabel{ppcurvature}
\RR (\del_i,\del_{+},\del_j,\del_{+})= \del_i\del_j H,
\end{equation}
and the components that are determined by this term via the symmetries of $\RR$. 
That is, we have
\[
\RR = 4\del_i\del_jH (\d x^i\wedge \d x^+ )(\d x^i\wedge \d x^+ ),
\]
in which we use Einstein's summation convention, and $\vf\wedge \psi=\frac{1}{2}(\vf\otimes \psi-\psi\otimes \vf)$  and $\vf  \psi=\frac{1}{2}(\vf\otimes \psi+\psi\otimes \vf)$ are the alternating and the symmetric product of two tensors.
Hence, the Ricci tensor of a pp-waves is given by
\[
\mathrm{Ric}= -\Delta H(\d x^+),
\]
where $\Delta=\sumi\del_i^2$ is the flat Laplacian.
Moreover, the covariant derivative of $\R$ is
\[
\nabla R= 4 \d H_{ij} (\d x^i\wedge \d x^+ )(\d x^i\wedge \d x^+ ),
\] including the differentials of the functions $H_{ij}:=\del_i\del_jH$. This shows that for a pp-wave to be a plane wave it requires $\del_i\del_j\del_k H=0$. Therefore, for a plane wave 
 the function $H$ is a quadratic polynomial in the $x^i$'s, i.e., in normal Brinkmann coordinates we have
\[
2 H(x^+,\x)=\x^\top S(x^+)\x
\]
where $\x$ denotes the column vector $(x^1, \ldots, x^n)$ and $S(x^+)$ is a symmetric $n\times n$-matrix depending on $x^+$. 
Plane waves satisfy the vacuum Einstein equations, i.e., are Ricci-flat if and only if $S(x^+)$ is traceless for all $x^+$.

A subclass of plane waves are the solvable Lorentzian symmetric spaces, called {\em Cahen-Wallach spaces} after \cite{cahen-wallach70}. As symmetric spaces, they satisfy $\nabla \R=0$ which forces the matrix $S$ to be constant. If $S$ has has no trace, Cahen-Wallach spaces provide remarkable examples of Ricci-flat, non flat symmetric spaces, contrasting the Riemannian situation where Ricci-flat symmetric spaces are flat.

The relation \eqref{ppcurvature} on a coordinate neighbourhood shows that the rank of $\R$ as an endomorphism  of $\Lambda^2T\M$ is equal to $n$ if and only if $\det(\hess(H))\neq0$. Indeed, the rank is smaller than $n$ if an only if there is a vector $X=\xi^i\del_i\in V^\perp$ such that $\R(X\wedge \del_+)=0$ which is equivalent to $\R(X,\del_+,\del_j,\del_+)=0$ for all $j$, i.e., $\xi^i\del_i\del_j H=0$.
 The curvature of a pp-wave and its derivatives are mapped into its holonomy algebra at $p$ as follows, where we work with normal Brinkmann coordinates centred at $p$:
\[
(\nabla_{X_1} \ldots \nabla_{X_k}\R)(\del_i,\del_+)\mapsto
\begin{pmatrix}
0& \left(X_1( \ldots (X_k(\del_i(\del_j H\ldots )|_0 \right)_{j=1}^n&0
\\
0&0&
\vdots
% -\left(X_1( \ldots (X_k(\del_i(\del_j(H))\ldots \right)_{j=1}^n
\\
0&0&0
\end{pmatrix}.
\]
This shows that if there is {\em one point} where the Hessian of $H$ has determinant not zero, then the pp-wave is indecomposable. However, the following example shows that the converse not true, i.e., there are indecomposable pp-waves, for which the rank of the curvature endomorphism is smaller than $n$ {\em on an open set}.
\bbsp\label{rankexample}
We give an example of a strongly indecomposable $4$-dimensional plane wave whose curvature has a kernel everywhere, and which therefore has  everywhere rank $1$. Given two functions $a_1$ and $a_2$ on $\rr$ with $a_1^2+a_2^2\neq0$ we consider the  matrix
\[
S=
\begin{pmatrix}
a_1 \\ a_2
\end{pmatrix}
\begin{pmatrix}
a_1 & a_2
\end{pmatrix}
=\begin{pmatrix}
a_1^2 &a_1a_2 \\
a_1a_2 & a_2^2
\end{pmatrix}.
\]
which has constant rank one.
Then $S$ defines a plane wave metric
\[
\gg=2\d x^+(\d x^-+ \x^\top S(x^+) \x \d x^+) +\d\x^2.
\]
Its curvature tensor is given by the matrix $S$ and hence has  everywhere rank $1$. However the derivative of the curvature is given by the matrix
\[
(\nabla_{\del_+}\R)(\del_+,\del_i,\del_+,\del_j)  =\dot a_i a_j +a_i\dot a_j
\]
which has determinant 
\[
\det(\dot S)= 4a_1a_2\dot a_1\dot a_2-\left(\dot a_1 a_2+a_1\dot a_2\right)^2=-\left(\dot a_1 a_2-a_1\dot a_2\right)^2,
\]
which in general is not zero. Therefore, as the first derivative of the curvature has no kernel, the holonomy of  $\gg$ is equal to $\rr^2$ and hence $\gg$ is strongly indecomposable.

We can even choose the matrix $S$ in a way that the resulting plane wave is homo\-geneous.
Indeed, if we
set
\[
S_- = \begin{pmatrix}
1 & 0 \\
0 & 0
\end{pmatrix},
\quad
F=\begin{pmatrix}
0 & -1 \\
1 & 0
\end{pmatrix},
\]
then
\[
\exp(x^+ F) =
\begin{pmatrix}
\cos(x^+) & -\sin(x^+) \\
\sin(x^+) & \cos(x^+)
\end{pmatrix}
\]
and
\[
S(x^+) = \exp(-x^+ F) S_- \exp(x^+ F)
=\begin{pmatrix}
\cos(x^+)^2 & -\cos(x^+)\sin(x^+)\\
-\cos(x^+)\sin(x^+)  & \sin(x^+)^2
\end{pmatrix},
\]
 has constant rank one.
According to Blau and O'Loughlin \cite{blau-oloughlin03}
(see also Section \ref{rem_hom_PW}),
$S$ defines a homogeneous plane wave metric and 
\[
\det(\dot S)=
-6\cos(x^+)^2\sin(x^+)^2
-\cos(x^+)^4-\sin(x^+)^4
\neq0
\]
shows that it is indecomposable.
\ebsp

In order to deduce Corollaries \ref{folg1} and \ref{folg3} from Theorem \ref{maintheo} we  observe 
\blem\label{ricrank}
If a pp-wave is Ricci-flat and its curvature endomorphism has
rank $1$, then it is flat.
\elem
\bprf
Assume that the curvature endomorphism has rank $1$ at a point $p$. This implies that there is an orthonormal
basis $V_p,E_1, \ldots , E_n$  of $V^\perp|_p$
such that $\R(E_+,E_i,E_+,E_j)=0$ unless $i=j=1$, where $E_+$ is transversal to $V^\perp|_p$. But then
\[
0=\mathrm{Ric} (E_+,E_+)=\sumi \R(E_+,E_i,E_+,E_i)=\R(E_+,E_1,E_+,E_1),
\]
so $\gg$ is flat.
\eprf
Regarding indecomposability, in what follows the following observation will be useful:

\blem\label{indeclemma}
On  a pp-wave $(\M,\gg)$, let $\cal U$ be simply connected patch of   Brinkmann coordinates and 
let $L=a^i\del_i$ be a non-zero vector field on $\M$ with constant coefficients $a^i$ such that $\RR(X,Y)L=0$ for all $X,Y\in T\cal U$. 
Then  the holonomy of $(\cal U,\gg)$ is properly contained in $\rr^n$, i.e., it does not act indecomposably. Moreover, $\gg$ is locally a product metric.
%, and, if we assume $\M$ to be simply connected and $\gg$ geodesically complete, $(\M,\gg)$ is globally a product.
\elem
\bprf Since $L=a^i\del_i$ has constant coefficients and no
$\del_+$-component, it is easy to see that its parallel transport  along a curve $\gamma$ is given as
$\mathrm P_\gamma(L|_{\gamma(0)})=\lambda \del_- +L|_{\gamma(1)}$ for some $\lambda\in \rr$ depending on the curve. 
Since $L$ as well as $\del_-$ are annihilated by the curvature tensor, we get that 
\[
\RR(X,Y)\circ \mathrm{ P}_\gamma (L) =\RR(X,Y)(\lambda \del_-+L)=0.\]
Using the Ambrose-Singer holonomy theorem, this shows that not only the null vector $\del_-$ but also the space-like vector $L$ is invariant under the holonomy algebra of $(\cal U,\gg|_{\cal U})$, which, as a consequence is reduced from $\rrn$ to a decomposable subalgebra. The reminder of the statement follows from the local version of the  de Rham--Wu decomposition theorem.\eprf
We conclude this section with a proof Theorem \ref{jek}. It generalises the proof in \cite{jek60} but avoids the use of coordinates.

\bprf[{\bf Proof of Theorem \ref{jek}}]
Let $(\M,\gg)$ be a semi-Riemannian manifold of dimension $n+2$ and $K_-,K_1, \ldots , K_n$ {\em commuting} Killing vector fields such that $K_-$ is null and the $K_i$ are orthogonal to $K_-$. We will show that this implies that $V:=K_-$ is parallel and that  $\R(X,Y)Z=0$ and $\nabla_X\R=0$ whenever $X,Y,Z\in V^\perp$. 
First note that we have
\be
\gg(\nabla_{K_i}K_j,K_k)
&=&
-
\gg(\nabla_{K_k}K_j,K_i)
=
-
\gg(\nabla_{K_j}K_k,K_i)
= 
\gg(\nabla_{K_i}K_k,K_j)
\\
&=&
\gg(\nabla_{K_k}K_i,K_j)
=
-
\gg(\nabla_{K_j}K_i,K_k)
=
-
\gg(\nabla_{K_i}K_j,K_k),
\ee
and hence \belabel{ijk}
\gg(\nabla_{K_i}K_j,K_k)=0\eeq for $i,j,k=0, \ldots , n$.
Set $\gg_{ij}:=\gg(K_i,K_j)$.
Clearly, $\gg_{0i}=0$ but the Koszul formula also gives
\belabel{dgijk}
\d\gg_{jk}(K_i)
=
\gg(\nabla_{K_i}K_j,K_k) + \gg(\nabla_{K_k}K_i,K_j)
=0.\eeq
Now we show that $V=K_-$ is parallel.
To this end fix a null vector field $Z$ such that $\gg(V,Z)=1$ and $\gg(Z,K_i)=0$ for $i=1, \ldots, n$. Clearly we have $\gg (\nabla_{Z}K_i,Z)=0$, and the Koszul formula also gives us that
\be
0
&=&
\gg (\nabla_{Z}K_i,K_j)+
\gg (\nabla_{K_j}K_i,Z)
\\
&=&
\frac{1}{2} 
\left(Z(\gg_{ij})-Z(\gg_{ij})+
\gg([Z,K_i],K_j)+\gg([K_J,Z],K_i) + \gg([Z,K_i],K_j)+\gg([Z,K_j],K_i)\right)
\\
& =& \gg([Z,K_i],K_j).
\ee
This implies that
\[
\gg (\nabla_{Z}K_i,K_j)
=
-
\gg (\nabla_{K_j}K_i,Z)
=\frac{1}{2}Z(\gg_{ij}),
\]
and in particular that $\nabla_ZV =0$ and $\nabla_{K_i}V=0$, i.e., that $V=K_-$ is parallel.
Moreover, we obtain that
\[\nabla_{K_i}K_j=-\frac{1}{2}Z(\gg_{ij}) V.
\]
This implies that
\belabel{ppeqn}
\begin{array}{rcl}
2 \R(K_i,K_j)K_k
&=&
\left(K_j(Z(\gg_{ik}))- K_i(Z(\gg_{jk})) \right)V\\
& = &
\left([K_j,Z](\gg_{ik})- [K_i,Z](\gg_{jk}) \right)V\\
& =& 0,
\end{array}
\eeq
because of \eqref{dgijk} and since the equation $0=\gg([Z,K_i],V)$ from above shows that $[Z,K_i]$ has no $Z$-component. Hence, we have shown that $\R(K_i,K_j)K_k=0$, i.e., that $\gg$ is a pp-wave in the case when $\gg$ is Lorentzian. In order to show that $\nabla_X\R=0$ for all $X\in V^\perp$ we use the integrability condition \eqref{kill-int}. Denote by $\phi_i:=\nabla K_i$. Obviously $\phi_-=0$ and $\phi_i(K_j)=-\frac{1}{2}Z(\gg_{ij}) V$ and $\phi_i(Z)\in \mathrm{span}(K_i)_{i=0}^n$. This and \eqref{ppeqn} together with  
the integrability condition \eqref{kill-int} gives us
\[\nabla_{K_i}\R=\phi_i\cdot \R=0,\]
and hence the statement of Theorem \ref{jek}.\eprf

%The following observation will be useful.
%\blem
%A pp-wave is strongly indecomposable if and only if the rank of the curvature considered as linear map  $\R : \Lambda^2T\M\to \Lambda^2T\M$ is constant and  equal to $n$ on a dense set in $\M$.
%\elem
%
%
%\bprf By the very definition of a pp-wave,  the kernel of $\R : \Lambda^2T\M\to \Lambda^2T\M$ has at least dimension $\frac{1}{2}n(n+1)+1$. hence its rank is at most $n$. Now assume there is an open $\cal U$ set where it is smaller than $n$. We introduce normal Brinkmann coordinates $(x^-,x^i,x^n)$ on this (or a smaller) set. If the rank of $\R$ is smaller than $n$ on this open set,  there is a vector field $X\in \Gamma(V^\perp|_\cal U)$  which is transversal to $V$ and such that $\R(X,\del_+)=0$ . siehe Lemma \ref{indeclemma} ...
%\eprf
\section{The Killing equation for pp-waves}
\label{killsec}
\subsection{The Killing equation in normal Brinkmann coordinates}
Here we derive the Killing equation in Brinkmann coordinates and then specialise this to normal Brinkmann coordinates found in Lemma \ref{coordlemma}. Mostly we follow \cite{blau-oloughlin03} where the Killing equation for plane waves in Brinkmann coordinates is derived and solved.
We fix Brinkmann cooridnates $(x^-, \x=(x^1, \ldots, x^n), x^+)$ and, using \eqref{nablai}, compute
the Lie derivative $\cal L_K\gg$ of the metric $\gg$ in direction of a vector field 
\[K:=K^-\del_-+K^i\del_i+K^+\del_+,\]
as
%let $\kappa^-$, $\kappa^i$ and $\kappa^+$ be one-forms defined by 
%\[\nabla K= \kappa^-\otimes \del_- + \kappa^i\otimes \del_i + \kappa^+\otimes\del_+.
%\]
%Using the $\kappa$'s we can express the bilinear form 
%\[
%\gg (\nabla K,\cdot )= \kappa^-\otimes \d x^+ + \delta_{ij}\kappa^i\otimes \d x^j +\kappa^+\otimes (\d x^-+2H \d x^+)
%\]
%and its symmetrisation 
%\[
%\frac{1}{2} \cal L_K\gg=
%\kappa^-   \d x^+ + \delta_{ij}\kappa^i  \d x^j +\kappa^+  (\d x^-+2H \d x^+).
%\]
%Now we compute, using \eqref{nablai} and introducing the notations 
%$H_i:=\del_i H$, $H^i:=\delta^{ij}H_j$, $\dot H:=\del_+ H$, and in general a dot for $\del_+$ derivatives,
%that
%\be
%\kappa^+&=&\d K^+
%\\
%\kappa^i&=& \d K^i-K^+ H^i \d x^+
%\\
%&=&\del_-K^i \d x^- +\del_jK^i \d x^j +(\dot{K}^i -K^+ H^i)\d x^+
%\\
%\kappa^-&=& \d K^- +K^iH_i \d x^+ +K^+ \d H\\
%&=&\del_-K^-\d x^-+ (\del_iK^-+ K^+H_i) \d x^i +(\dot{K}^-+ K^iH_i+K^+ \dot H)\d x^+.
%\ee
%This implies
\be
\tfrac{1}{2}\cal L_K\gg
&=& 
\del_-K^+(\d x^-)^2 
+
\delta_{ij}\del_k K^i\d x^k\d x^j
+
\left( \dot{K}^-+K^iH_i +K^+\dot H +2 H\dot{K}^+\right) (\d x^+)^2 
\\
&&+
\left(\delta_{ij}\del_-K^j+\del_iK^+\right) \d x^-\d x^i
+
\left(\del_iK^-+\dot{K}^i+2 H\del_iK^+\right)\d x^i\d x^+
\\
&&
+
\left( \del_-K^-+ 2H  \del_-K^+ +\dot{K}^+\right) \d x^-\d x^+,
\ee
where we write $H_i:=\del_i H$, $H^i:=\delta^{ij}H_j$, $\dot H:=\del_+ H$, and in general a dot for $\del_+$ derivatives.
Hence, $K$ is a Killing vector field if and only if its components satisfy the following system
\begin{eqnarray}
\del_-K^+&=&0 \label{--}
\\
\del_iK^j+\del_jK^i&=&0\label{ij}
\\
\dot{K}^-+K^iH_i +K^+\dot H +2 H\dot{K}^+&=&0\label{++}
\\
\del_-K^i+\del_iK^+&=&0\label{-i}
\\
\del_iK^-+\dot{K}^i+2 H\del_iK^+&=&0
\label{+i}
\\
\del_-K^-+\dot{K}^+&=&0
\label{-+}
\end{eqnarray}
These equations were derived  in \cite{blau-oloughlin03} and in the following we  review some of the arguments given there. Because of \eqref{--}, $K^+$ is independent of $x^-$. 
Hence, 
 differentiating 
\eqref{-i} and \eqref{-+} with respect to $x^-$ gives
\[
\del_-^2K^i=\del_-^2K^-=0,\]
showing that $K^-$ and all $K^i$ are linear in $x^-$, whereas differentiating \eqref{ij} with respect to $x^-$ and 
\eqref{-i} with respect to $x^j$ and symmetrising over $i$ and $j$ gives
\[
0=
%\del_j\del_-K^i+\del_i\del_-K^j+2\del_i\del_jK^+
%=
2\del_i\del_jK^+
\] for all $i,j$ showing that $K^+ $ is linear in the $x^i$'s. Hence, there are functions $\alpha^+,\alpha_1 \ldots , \alpha_n$ of $x^+$ only such that
\[
K^+=\alpha_ix^i +\alpha^+.
\]
Then equation \eqref{-+} becomes
\[
0=
\del_-K^-+\del_+K^+ =\del_-K^-+ \dot \alpha_i x^i + \dot\alpha^+,
\]
and hence, there is a function $A^-:=A^-(x^+,x^1, \ldots, x^n)$ depending on $(x^+,x^1, \ldots, x^n)$ such that
\[
K^-=-(\dot \alpha_i x^i + \dot\alpha^+)x^-+A^-.
\]
Furthermore, equation \eqref{-i} becomes
\[
0=
\del_-K^i + \alpha_i
\]
yielding the existence of functions $A^i=A^i(x^+,x^1, \ldots, x^n)$ depending on $(x^+,x^1, \ldots, x^n)$ such that
\[
K^i= -\alpha_i x^- +A^i.
\]
With this information at hand, we evaluate 
\eqref{+i} and get
\[
0
=
-2\dot\alpha_i x^-+\del_iA^- +\dot A^i+2 H \alpha_i.
\]
Since $A^-$, $A^i$ and $H$ are independent of $x^-$ this shows that the $\alpha_i$'s are constant, i.e., $\alpha_i\equiv a_i\in \rr$ . Hence, $K$ is a Killing vector field if and only if its components are given as
\be
%\label{k+}
K^+&=&a_ix^i +\alpha^+
\\
%\label{k-}
K^-&=&-\dot\alpha^+ x^-+A^-
\\
%\label{ki}
K^i&=& -a_i x^- +A^i
\ee
for constants $a_i$, a function $\alpha^+$ of $x^+$ and functions $A^-$ and $A^i$ of $(x^1, \ldots , x^n,x^+)$ subject to the equations
\begin{eqnarray}
\label{++A}
-( \ddot \alpha^+ +a_iH^i)x^- +\dot A^- +A^iH_i +(a_ix^i+\alpha^+)\dot H +2 H\dot \alpha^+
&=&0
\\
\label{ijA}
\del_iA^j+\del_jA^i&=&0
\\
\del_iA^- +\dot A^i+2 H a_i&=&0.
\label{+iA}
\end{eqnarray}
Differentiating \eqref{++A} with respect to $x^-$ and then with respect to $x^i$ we obtain
\belabel{kerHess}
a_i\del_j\del^i H=0.
\end{equation}
Recalling formula \eqref{ppcurvature}, this shows that the vector field
$L=a^i\del_i$ on $\M$, for $a^i:=a_i$ constants,  is annihilated by the curvature tensor $\RR$ of $\gg$, i.e., $\RR(X,Y)L=0$ for al $X,Y\in T\M$. 

From now on we will  assume that $(\M,\gg)$ is strongly indecomposable, i.e., that the   holonomy algebra of $(\cal U,\gg|_{\cal U})$ acts indecomposably.  Under this assumption, Lemma \ref{indeclemma}  implies by \eqref{kerHess} that all the constants $a_i$ vanish,
\[
a_i=0.
\]
Differentiating equation \eqref{++A} with respect to $x^-$ yields
that $\alpha^+=ax^++b$ is linear.

Now, with the $a_i$ being zero, differentiating 
equation \eqref{ijA} with respect to $x^+$ and equation 
\eqref{+iA} with respect to $x^j$ and symmetrising over $i$ and $j$ gives us
\[\del_i\del_iA^-=0,\]
which shows that $A^-$ is linear in the $x^i$'s.
 Plugging this back into \eqref{+iA}, and differentiating with respect to $x^j$ yields
 
\[
0= \del_j\dot A^i,\]
which shows that $A^i$ is of the form
$A^i=\psi^i+F^i$, where the 
 $\psi^i$ are functions of $x^+$ only and $F^i$ are functions of $(x^1, \ldots , x^n)$.
 Consequently, there is a function $\vf$ of $x^+$ such that
 \[
 A^-(x^+,\x)= -\x^\top \dot\Psi+\vf.
 \]
 where we write $\Psi=(\psi^1, \ldots , \psi^n)$.
 
 Finally, the functions $F^i$ are subject to the Euclidean Killing equation
\[
\del_i F^j+\del_j F^i=0,
\]
 the solutions of which are given, up to constants, by a skew-symmetric matrix $f^i_{~j}=-f^j_{~i}$ such that
 $F^i=f^i_{~j}x^j$. 
 Plugging all this back into equation \eqref{++} we 
 obtain that
 any Killing vector field $K$ on an  indecomposable pp-wave $(\M,\gg)$ in Brinkmann coordinates is of the form
\belabel{killing1}
K(x^-,x^+,\x)=-\left(a x^ - + \vf (x^+) +  \x^\top\dot\Psi(x^+)  \right)\del_- + \left(\Psi(x^+) +F \x \right)^i \del_i +(ax^++b) \del_+,
\end{equation}
 where $a$, $b$ and $F=(f^i_{~j})\in\so(n)$ are constant, and $\vf$ and $\Psi=(\psi^1, \ldots , \psi^n)$ are functions of $x^+$ satisfying the equation
 \belabel{killingeqn1}
- \ddot\Psi^\top \x-\dot \vf+\grad( H)^\top(\Psi +F \x)  + (ax^++b)\dot H+ 2a H=0
\end{equation}
Now, in normal Brinkmann coordinates, we can simplify equation \eqref{killingeqn1}:
 \btheo \label{killtheo}
 Let $(\M, \gg)$ be a strongly indecomposable  pp-wave,  $p\in \M$, and let   $(\cal U, (x^+,x^-,\x=(x^1, \ldots , x^n)))$ be  normal Brinkmann coordinates centred  at $p$ with $2H:=\gg(\del_+,\del_+)$.
Then $K$ is a Killing vector field if and only if 
\begin{equation}
K=(c- ax^- -  \dot \Psi^\top\x)\del_- + \left(\Psi+F \x\right)^i \del_i +(ax^++b) \del_+,
\label{killing}
\end{equation}
where $a,b,c\in\mathds{R}$, $F\in\so(n)$ are constant
and $\Psi\in C^\infty(\mathds{R},\mathds{R}^n)$ subject to the 
  Killing equation
\belabel{killingeqn}
\ddot\Psi^\top\x-\grad(H)^\top(\Psi+F \x) - (ax^++b)\dot H-2a H=0.
\end{equation}
Moreover, for the commutator $\hat{K}=[K_1,K_2]$ of two Killing vector
fields $K_1,K_2$ the parameters are
\begin{equation}
\begin{split}
\hat{a} &= 0\\
\hat{b} &= a_2 b_1-a_1 b_2\\
\hat{c} &=
\dot{\Psi}_1^\top \Psi_2- \Psi_1^\top\dot{\Psi}_2-a_1c_2+a_2 c_1\\
\hat{F} &= -[F_1,F_2]\\
\hat{\Psi} &=
F_2\cdot\Psi_1-F_1\cdot\Psi_2 +(a_1x^++b_1)\dot{\Psi}_2
-(a_2 x^+ +b_2)\dot\Psi_1.
\end{split}
\label{bracket}
\end{equation}
\etheo

\bprf
Clearly, $K$ in \eqref{killing} is a Killing vector field as its components satisfy equation \eqref{killingeqn1} with $\vf(x^+)\equiv -c$.

On the other hand, we have seen that every Killing vector field in Brinkmann coordinates is of the form \eqref{killing1} with components satisfying equation \eqref{killingeqn1}. Choosing the Brinkmann coordinates to be normal at $p$, equation 
 \eqref{killingeqn1} when taken along $\x=\oo$ becomes $\dot\vf\equiv 0$, which we solve by $\vf(x^+)\equiv -c$.

Finally, it is a matter of checking that the induced Lie bracket is of the form \eqref{bracket}.
Note that, as required,  the term  $\dot{\Psi}_1^\top\Psi_2-\Psi_1^\top \dot{\Psi}_2$ is constant as a consequence of  both $\Psi_1$ and $\Psi_2$ being solutions of equation \eqref{dkillingeqn0}.
\eprf
Let us make a few observations. The fact that $c$ does not appear in \eqref{killingeqn} is due to  $\del_-$, as a parallel vector field,  is a Killing vector field. 
Moreover, the parameters $a,b,c,F,\Psi$ uniquely determine the Killing vector field $K$, which is  determined by the values of its covariant derivative at the point $p$.
For the covariant derivatives of $K$
we compute
\belabel{nabK}
\begin{array}{rcl}
\nabla_{\del_-}K&=& -a \del_-
\\
\nabla_{\del_i} K&=&-\left( \dot \psi^i-(ax^++b)\del_iH\right)\del_- +f_{i}^{~k}\del_k
\\
\nabla_{E_+}K&=& \left( \dot \psi^i-(ax^++b)\del_iH\right)\del_i + aE_+,
\end{array}\eeq
where $E_+=\del_+-H\del_-$ and we have to use the Killing equation \eqref{killingeqn1} to obtain the last derivative. Hence, 
at zero, the Killing vector in \eqref{killing1} and its covariant derivative is given by 
\belabel{killinit}
\begin{array}{rcl}
K|_0&=&c \del_-+\psi^i(0)\del_i +b\del_+
\\
\nabla_{\del_-}K|_0&=&-a\del_- 
\\
\nabla_{\del_i}K|_0&=&-\dot \psi^i(0)\del_- +f_{i}^{~k}\del_k
\\
\nabla_{\del_+}K|_0&=&
\dot\psi^i(0)\del_i +a\del_+
\end{array}
\eeq
Moreover, differentiating equation \eqref{killingeqn} yields
\beq
\label{dkillingeqn}
\ddot\Psi+  F \grad( H)-  \hess(H) (\Psi+F \x) -(ax^++b)\grad(\dot H)-2a \grad( H)=0.
\eeq
By the properties of the normal Brinkmann coordinates from Lemma \ref{coordlemma}, this becomes a second order linear ODE system for
$\Psi=(\psi^1, \ldots , \psi^n)$
when taken along $\x=\oo$:
\belabel{dkillingeqn0}
\ddot{\Psi}(t)
-\hess(H)(t,\oo) \Psi (t)
=0.
\end{equation}
Fixing initial conditions $\Psi(0)$ and $\dot{\Psi}(0)$ gives  a unique solution to this system. This illustrates how $K$ is completely determined by the initial conditions.

In the remainder of the section we will  consider some special cases, known results and examples.

\subsection{Transversal Killing vector fields} 
\label{transverseremark}
We will see that a crucial issue of the Killing equation on pp-waves is the existence of Killing vector fields that are transversal to the parallel null distribution $V^\perp$ of
rank $n+1$.

First note that, if  
$\dot H=0$,  then there is always the transversal Killing vector field $\del_+$, but in general transversal Killing vector fields are much harder to find and the situation is much more involved. For example, for certain pp-waves there exist Killing vector fields with 
 $b=0$ but $a\neq0$ being  tangent to $V^\perp$ {\em only} along the leaf $x^+=0$ but  transversal elsewhere, i.e., pp-waves for which 
\[\k(V):=\{ K\in \k\mid \gg(K,V)=0\}\]
and 
\[\k_p(V):=\{ K\in \k\mid \gg(K,V)|_p=0\}\]
are different.
Note that Theorem \ref{killtheo} and formulae \eqref{nabK} show that 
\[\k':=\{ K\in \k\mid \nabla_VK=0\}\]
and its subalgebra $\k(V)$
are actually  ideals in the Lie algebra $\k$ of Killing vector fields. In fact we have that $[\k,\k]=\k'$.
Killing vector fields that are transversal {\em at some point} project onto non-zero elements in 
 the quotient
Lie algebra $\k/\k(V)$.
\bfolg\label{transcol}
The Lie algebra $\k/\k(V)$ is isomorphic to a subalgebra of $\mathfrak{aff}(1)$, the Lie algebra of affine transformations  of $\rr$. In particular, if $\k/\k(V)$ is $2$-dimensional, then there are two Killing vector fields $K$ and $\hat K$ such that
\be
K &=& x^+\del_+ \mod V^\perp
\\
\hat K &=& \del_+ \mod V^\perp.
\ee
\efolg
\bprf
The theorem shows that there is a Lie algebra homomomorphism
\[\k \ni (ax^++b)\del_++ K^i\del_i+K^-\del_-\mapsto (a,b)\in \mathfrak{aff}(1),\]
the kernel of which is $\k(V)$. Hence $\k/\k(V)$ injects  homomorphically into $\mathfrak{aff}(1)$. If $\k/\k(V)$ is $2$-dimensional, we can invert this map obtaining  two Killing vector fields of the required form.
\eprf

\begin{bsp} 
\label{ex_dim3}
Here we will give an example of a  $3$-dimensional pp-wave for which the Lie algebra $\k/\k(V)$ is indeed $2$-dimensional, and more importantly, which is {\em locally homogeneous but  not a plane wave}, showing that the assumption on the curvature in Theorem~\ref{maintheo} is essential. %and, in particular, that Theorem \ref{maintheo} does not cover pp-waves of dimension $3$.
Consider the  pp-wave $(\M,\gg)$ where $\M=\mathds{R}^3$
and
\[
\gg = 2\d x^+(\d x^- + \mathrm{e}^{2a x} \d x^+) + \d x^2,
\]
where $a\in\mathds{R}\setminus\{0\}$ is a constant and
$(x^+,x^-,x)$ are the standard coordinates in $\mathds{R}^3$.
In particular, the function $H(x^+,x)$ is 
\[
H(x) = \mathrm{e}^{2a x}.
\]
Since  $\frac{\partial H}{\partial x^+}=0$, the
Killing equation \eqref{killingeqn} takes the form
\[
\ddot{\psi}x - 2a\mathrm{e}^{2a x}\psi
-2a\mathrm{e}^{2a x} = 0.
\]
Solving this equation, we find that, in addition to $V=\del_-$ and $\del_+$, there is another Killing vector field of $\gg$, namely
\begin{eqnarray*}
%K_-&=&V\ =\ \partial_-, \\
%K_+&=&\partial_+, \\
K&=&a x^+\partial_{+} - a x^- \partial_{-} - \partial_x.
\end{eqnarray*}
Hence, $\k$ is $3$-dimensional. Since $\gg(K,V)=a x^+$, we have $\k(V)=\rr\cdot \del_-$ and thus $\dim(\k/\k(V))=2$. Moreover, 
the Killing vector fields vector fields span the tangent space $T_p \M$ at any
point $p\in\M$, so $(\M,\gg)$ is a locally homogeneous
pp-wave.
However,  $(\M,\gg)$ is strongly indecomposable since 
\[\R(\del_x,\del_+)=
\begin{pmatrix}
0&2a \mathrm{e}^{2a x} &0
\\
0&0&
-2a \mathrm{e}^{2a x} 
% -\left(X_1( \ldots (X_k(\del_i(\del_j(H))\ldots \right)_{j=1}^n
\\
0&0&0
\end{pmatrix}
\not= 0,
\]
for any $x\in \rr$, and
$(\M,\gg)$ clearly is {\em not a plane wave} since
\[(\nabla_{\del_x}\R)(\del_x,\del_+)=
\begin{pmatrix}
0&4a^2 \mathrm{e}^{2a x} &0
\\
0&0&
-4a^2 \mathrm{e}^{2a x} 
% -\left(X_1( \ldots (X_k(\del_i(\del_j(H))\ldots \right)_{j=1}^n
\\
0&0&0
\end{pmatrix}
\not= 0.
\]
\end{bsp}

\bbsp[Ehlers \& Kundt]\label{transex}
Similar examples with $\dim( \k/\k(V) )=2 $ but in dimension $4$ are given by Ehlers and Kundt
in \cite[Table 2-5.1]{ehlers-kundt62} as a correction to \cite{jek60}. For one class of examples $H$ is given as the real part of the complex function
\[
\mathrm{e}^{2a z}, \quad \text{ with }a>0,
\]
of $z=x^1+\mathrm{i}x^2$.
Then $\del_-$ and $\del_+$ and 
\[-a(x^-\del_-+x^+\del_+)-\del_1\]
span the Killing vector fields.
For the other class, $H$ is given as the real part of 
\[
\mathrm{e}^{2\mathrm{i} a \ln(z)}, \quad \text{ with }a\not=0.
\]
Here, the Killing vector fields are spanned
by $\del_-$ and $\del_+$ and 
\[-a(x^-\del_-+x^+\del_+)+x^1\del_2-x^2\del_1.\] 
Note that with $\dim(\k)=3$ and $\dim(\k_p(V))=2$ both metrics are neither homogeneous nor $V^\perp$-homogeneous.
\ebsp

\bbsp[Sippel \& Goenner]\label{sgex}
Another example of this type with  $\dim( \k/\k(V) )=2 $ in dimension $4$ was given by
Sippel and Goenner
in \cite[Table II, no. 9]{sippel-goenner86}. These examples are  pp-wave metrics on $\rr^4$ which are locally homogeneous but not plane waves. However, they turn out to be {\em decomposable}. 
The pp-wave metric is defined by
\[H(x^1,x^2):=c\,\mathrm{e}^{a_1x^1-a_2x^2},\]
with $c,a_1,a_2$ constants with $a_1^2+a_2^2\neq0$.
The Killing vector fields are given by $\del_-$, $\del_+$ and 
\be
K&:=&x^+(a_2\del_1+a_1\del_2)+ (a_2x^1+a_1x^2)\del_-\ \in\  \k(\del_-),\\
K_i&:=&\del_i+a_i(x^+\del_+-x^-\del_-),
\ee
for $i=1,2$,
and span the tangent space. However a coordinate transformation \[x=a_1x^1-a_2x^2,\ \  y=a_2x^1+a_1x^2\] reveals that this metric is decomposable.
\ebsp

For plane waves  we can show
\bs\label{transprop}
Let $(\M,\gg)$ be a strongly indecomposable   plane wave. Then 
\belabel{kmodkv}
\dim (\k/\k(V))\le 1.\eeq
\es
\bprf
Assume there are two linearly independent Killing vector fields that are not tangent to $V^\perp$.  They are of the form
\be
K&=& x^+\del_+ + (\psi +F \x)^k\del_k + K^-\del_-
\\
\hat K&=& \del_+ + (\hat\psi +\hat{F} \x)^k\del_k + \hat{K}^-\del_-.
\ee
Now, differentiating equation \eqref{dkillingeqn} again we obtain
\belabel{ddkillingeqn}
(\psi +F \x)^k\del_k\hess(H) + [F,\hess(H)] 
+
(ax^++b) \hess(\dot H)
+2a \hess(H)=0.
\eeq
For a plane wave in normal Brinkmann coordinates with $S=\hess(H)$ we have that $\del_kS=0$ and thus when taking  equation \eqref{ddkillingeqn} along $\x=\oo$, we obtain for $K$ and $\hat K$ that
\begin{align*}
[F,S]-x^+ \dot{S}-2S &=0 \\ 
[\hat{F},S] -\dot S &=0.
\end{align*}
This implies that
\[
[F-x^+\hat{F} ,S]-2S=0,
\]
for all $x^+$. 
Since the map $S\mapsto[F-x^+\hat{F},S]$ when acting on symmetric matrices
 is skew-symmetric  with respect to the trace form, which, on the other hand, is positive definite on symmetric matrices, we obtain that $S\equiv 0$, which is a contradiction.
 \eprf
 A fundamental question is whether, in dimensions greater than $3$,  ($V^\perp$-) homogeneity and inecomposability forces $\k/\k(V)$ to have dimension $1$. Because of the additional term $\del_k\hess(H) $, we are not able to prove \eqref{kmodkv} for arbitrary ($V^\perp$-) homogeneous pp-waves, but we conjecture that it is true:
\begin{conj}\label{trans-con}
For an indecomposable locally homogeneous pp-wave of dimension greater than $3$, the Lie algebra $\k/\k(V)$ is at $1$-dimensional.
\end{conj}
Our proof of Theorem \ref{maintheo} in Section \ref{proofsec} will show that if this conjecture is true, then in dimensions greater than $3$ we can drop the assumption on the curvature in Corollary \ref{folg2} (see Remark \ref{con-rem}).

% Because of the additional term $\del_k\hess(H) $, we are not able to prove \eqref{kmodkv} for arbitrary ($V^\perp$-)homogeneous pp-waves but believe that it is true:
%%\begin{conj}\label{trans-con}
%%For an indecomposable locally homogeneous pp-wave $\k/\k(V)$ is at most $1$-dimensional.
%%\end{conj}
%%Our proof of Theorem \ref{maintheo} in Section \ref{proofsec} will show that if this conjecture is true, we can drop the assumption on the curvature in Corollary \ref{folg2} (see Remark \ref{con-rem}).

\subsection{Plane waves}
\label{pw-sec}
 In this section we will recall some facts about plane waves, for which the  Killing equation is completely solved in \cite{blau-oloughlin03}.
\subsubsection{Plane waves and the Heisenberg algebra}\label{rem_planewave_heisenberg}
For a plane wave defined by a matrix $S(x^+)$ the Lie algebra 
 $\k(V)$ always contains the Heisenberg algebra  $\he(n)$. Indeed, for a plane wave we have 
 \[
 H=\tfrac{1}{2}\,\x^\top S(x^+) \x\] 
for a symmetric $x^+$-dependent matrix $S$, and hence
 \[\grad(H)=S \x\ , \ \ \mathrm{Hess}(H)=S.\]  
For such $H$, multiplying the differentiated equation \eqref{dkillingeqn} by $\x$ 
 implies the   Killing equation \eqref{killingeqn}, which therefore becomes equivalent to
  \eqref{dkillingeqn}. On the other hand, when setting $F=0$ and $a=b=0$, equation   \eqref{dkillingeqn} is equivalent to the linear ODE system 
 \eqref{dkillingeqn0} which, for a plane wave, becomes
 \belabel{dkillingeqn0pp}
 \ddot\Psi- S \Psi=0.
 \eeq
 Hence, 
  we have Killing vector fields
\belabel{heisenberg}
\begin{array}{rcl}
L_i&:=& \phi^k_i\del_k -\x^\top\dot\Phi_i \cdot \del_- 
\\
K_i&:=& \psi^k_i\del_k-\x^\top \dot\Psi_i \cdot \del_- ,
\end{array}
\eeq
where
$\Phi_i=(\phi_i^k)_{k=1,\ldots ,n}$ and $\Psi_i=(\psi_i^k)_{k=1,\ldots ,n}$ are solutions to 
the linear ODE system \eqref{dkillingeqn0pp}
with initial conditions
\be
\Phi_i(0)=\oo,&& \dot \Phi_i(0)=\e_i
\\
\Psi_i(0)=\e_i,&& \dot \Psi_i(0)=\oo,
\ee
which span $\he(n)$.
Clearly, $\del_-$ commutes with the $K_i$'s and $L_j$'s and we have
\belabel{klbracket}
 [ L_i,K_j] =( \Phi_i^\top\dot\Psi_j - \Psi_j^\top\dot\Phi_i )\del_-=-\delta_{ij}\del_-
 \eeq
because the term $\Phi_i^\top\dot\Psi_j - \Psi_j^\top\dot\Phi_i$ is constant as a consequence of equation \eqref{dkillingeqn0}.

Clearly, 
for  plane waves,  there are {\em commuting} Killing vector fields $X_1,\ldots, X_n,\del_-$ spanning the null distribution $V^\perp$. Theorem \ref{jek} shows that this can {\em only} happen for plane waves.

\subsubsection{Homogeneous plane waves}\label{rem_hom_PW}
For plane waves, the Killing equation \eqref{dkillingeqn}
becomes the following matrix ODE:
\belabel{planedkillingeqn1}
[ S(x^+), F]  + (ax^++b)\dot S(x^+)+ 2a\, S(x^+) =0.
\eeq
In  Section \ref{rem_planewave_heisenberg} we saw that $\k$
always contains a Heisenberg algebra.
Now, for a plane wave to be locally homogeneous, we need an additional Killing vector field $K$ transversal to $V^\perp|_p$.
Hence, when working with normal Brinkmann coordinates centred at $p$, one has to find a solution of  equation \eqref{planedkillingeqn1} with
$b\not= 0$.
This was done  by Blau and
O'Loughlin in \cite{blau-oloughlin03}.
Depending on $a$ being zero or not, 
they found two families of homogeneous plane waves, where the
metrics in both families are determined by the choice of a
constant symmetric matrix $S_-$ and a constant skew-symmetric
matrix $F$.

In the first case, when $a=0$ we can assume $b=1$ and hence the Killing equation \eqref{planedkillingeqn1} just becomes
\[
[ S(x^+), F]  + \dot S(x^+)=0.
\]
Clearly this is solved by 
\[
S(x^+)= \mathrm{e}^{x^+ F}S_-\mathrm{e}^{-x^+ F}
\]
with a constant skew symmetric matrix $F$ and a constant symmetric matrix $S_-$.
Hence, the metrics in the first family are of the form
\belabel{type1}
\gg = 2\d x^+ \d x^- + (\x^\top\mathrm{e}^{x^+ F}S_-\mathrm{e}^{-x^+ F}\x) (\d x^+)^2 + \d\x^2.
\eeq
When defined on all of  $\rr^{n+2}$ they are geodesically complete (see for example 
results by Candela, Flores and S\'{a}nchez \cite[Prop. 3.5]{candela-flores-sanchez03}).

In the second case we have $a\not=0$ so that we can assume $a=1$.
Here 
 the Killing equation \eqref{planedkillingeqn1}  becomes an ODE with singularity at $x^+=-b$, 
\[ (x^++b)\dot S(x^+)+[ S(x^+), F]  + 2\, S(x^+) =0.\]
It has the solution
\[
S(x^+)
=
\frac{1}{(x^++b)^2}(\mathrm{e}^{\log(x^++b) F}S_-\mathrm{e}^{\log(-(x^++b)) F}),
\]
again for constant (skew) symmetric matrices $F$ and $S_-$.
Hence, homogeneous plane wave  metrics in the second family are of the form
\belabel{type2}
\gg = 2\d x^+ \d x^- + \frac{1}{(x^++b)^2}(\x^\top\mathrm{e}^{\log(x^++b) F}S_-\mathrm{e}^{\log(-(x^++b)) F}\x) (\d x^+)^2 + \d\x^2,
\eeq
for constants $F$, $S_-$ and $b$. They are only defined for $x^+>-b$ and hence
 geodesically incomplete. Clearly, metrics for different $b$  can be pulled back  by a translation $x^+\mapsto x^++b$ to the metric with $b=0$  on $\{x^+>0\}$. Hence, metrics with different $b$  are  isometric to each other.

\subsubsection{Reductivity of homogeneous plane waves}
\label{pwred-sec}
Here we will show that homogeneous plane waves are always \emph{reductive}.
This means that for some subalgebra $\k_0$ of $\k$ generating
a (locally) transitive group action, the stabiliser
$\h:=\{K\in \k_0\mid K|_p=0\}$ in $\k_0$ of a point $p$
has a vector space complement $\m$ in $\k_0$
with $[\h,\m]\subseteq\m$.
\bs
Homogeneous plane waves are reductively homogeneous.
\es
\bprf For a homogeneous plane wave, we take $\k_0$ to be the subalgebra
generated by the Killing fields
\[
K_+, \partial_-, K_1,\ldots,K_n, L_1,\ldots,L_n,
\]
where
$K_i$, $L_j$ are defined in \eqref{heisenberg} and 
$K_+=-ax^- \del_- +(F \x)^i \del_i + (ax^++b) \del_+$
for a certain $F=(f_{~i}^j)\in\so(n)$ is transversal to $V^\perp$, which exists for homogeneous plane waves
according to \cite[(2.42)]{blau-oloughlin03}.
Working at $p$ with normal Brinkmann coordinates centred at $p$, we see that $\h$ is spanned by
the $L_i$'s defined in  \eqref{heisenberg}. Then the $\h$-invariant complement $\m$ is spanned by $\del_-$, $K_+$ and the $n$ Killing vector fields
\[
M_i:=[K_+,L_i].
\]
Note that this implies that 
\[M_i|_p=b\dot \phi^k_i(0)\del_k|_p=b\del_i|_p\] 
Hence, since also $K_+|_p=b\del_+|_p$, the vector space $\m$ defined in this way is indeed a complement to $\h$. Moreover, since both $M_i$ and $L_j$ are tangent to $V^\perp$ and without rotational component we obtain from \eqref{bracket} that
\[
[L_j,M_i]= c \del_-
\]
for a constant $c$. 
 Therefore we have
$
[\h,\m]\subseteq\m
$
and the plane wave is  reductive.
\eprf
 
\subsubsection{Cahen-Wallach spaces} For Cahen-Wallach spaces, the matric $S(x^+)$ is constant and thus equation \eqref{planedkillingeqn1} can always be solved by setting $F=0$ and $a=0$ and thus yielding a Killing vector field transversal to $V^\perp$. For Cahen-Wallach spaces generically  the algebra of Killing vector fields contains the oszillator algebra $\rr\ltimes \he(n)$ and hence has dimension at least $2n+2$. The stabiliser algebra is equal to the holonomy algebra which is $\rrn$. A Cahen-Wallach space may have additional Killing vector fields in addition to $\rr\ltimes \he(n)$. In fact, the additional Killing vector fields are isomorphic to the centraliser in $\so(n)$ of the constant matrix $S$ defining the Cahen-Wallach space. Hence, it might have at most $\frac{1}{2}n(n-1)$ additional symmetries.

% at least in dimension four. Their proof, however generalises to any dimension. \marg{That's what we believe ...}
%\bs[{\cite[Theorem 4.5.2]{jek60}}]\label{jek}
%Let $(\M,\gg)$ be a Lorentzian manifold of dimension $n+2$. Assume there exists {\em commuting} Killing vector fields $V,X_1, \ldots ,X_n$ everywhere spanning $V^\perp$. Then $(\M,\gg)$ is a plane wave.
%\es

\subsection{Dimension four}
\label{dim4rem} In \cite{jek60} the Killing equation \eqref{killingeqn} for $4$-dimensional pp-waves is explicitly solved   under the assumption that $(\M,\gg)$ is Ricci-flat, i.e., that $H$ is harmonic, so that methods from complex analysis can be used. In particular, 
in \cite[table on p.~79]{jek60}, the dimension of the space of Killing vector fields of a $4$-dimensional,  indecomposable, {\em Ricci-flat} pp-wave have been determined  as $\dim(\k)= 1, \ 2,\ 3,\ 5,\ 6$, and 
 the metrics are explicitly given for each case.
 Moreover, in \cite{sippel-goenner86} the assumption of Ricci-flatness was dropped and  new algebras of dimension $5$, $6$ and $7$ appeared, almost reaching the upper bound of  $8$ we will deduce from Theorem~\ref{killtheo} in Corollary \ref{dim-cor}. 
  Further results about symmetries of $4$-dimensional pp-waves were obtained in \cite{AichelburgBalasin96,AichelburgBalasin97}.

\section{Proof of the main results}
\label{proofsec}
In this section we will draw the conclusions from Theorem \ref{killtheo} that eventually will lead to a proof of Theorem \ref{maintheo}. We assume that $(\M,\gg)$ is an indecomposable pp-wave with parallel null vector field $V$.
First we note:
\bfolg\label{nabvk-cor}
Each Killing vector field satisfies 
$\nabla_VK\in \rr V$.
\efolg
\bprf
Let $p\in \M$ be an arbitrary point and chose normal Brinkmann coordinates centred at $p$ (Lemma \ref{coordlemma}).
In these coordinates a Killing vector field $K$ is of the form \eqref{killing} with its covariant derivative given in \eqref{nabK}. Since $V=\del_-$ on the coordinate patch, we get $\nabla_VK=a\cdot V$.
\eprf

Now denote by $\k$ the Killing vector fields of $(\M,\gg)$. 
We describe the evaluation map $\kappa$ at a point $p\in \M$ at which we choose a basis 
\[
(E_-, E_1, \dots, E_i,E_+)
\]
of $T_p\M$ such that 
\[
\gg_p(E_-, E_+)=1,\quad \gg_p(E_i,E_j)=\delta_{ij},
\]
where $i,j=1, \ldots, n$, and all other $\gg_p(E_\alpha,E_\beta)=0$ for $\alpha,\beta\in \{-,+,1,\ldots, n\}$. 
Moreover, in the proofs we will use 
 normal Brinkmann coordinates centred at $p$ and such that  
\beq\label{basis}E_-=\del_-|_p,\quad
E_i=\del_i|_p,\quad
E_+=(\del_+-H\del_-)|_p=\del_+|_p.\eeq
In Theorem \ref{killtheo} we have seen that, for a Killing vector field $K$ there are real numbers $a,b,c,X^i,Y^i,F=(f_i^{~j})\in \so(n)$ such that
\belabel{atp}
\begin{array}{rcl}
K|_p&=&cE_-+X^iE_i+bE_+
\\
\nabla_{E_-}K|_p&=&-a E_-\\
\nabla_{E_i}K|_p&=&-Y_iE_-+f_{i}^{~k}E_k
\\
\nabla_{E_+}K|_p&=& Y^i E_i +aE_+.
\end{array}
\eeq
Furthermore, we write $Y=(Y_i)$, $X^\top=(X_i)$ for the row vectors and $X=(X^i)$, $Y^\top=(Y^i)$ for the column vectors.

If we denote by $v\in \rr^{1,n+1}$ the null vector that is the image of $V$ under the evaluation map $\kappa$, i.e. $\kappa(V)=(0,v)\in \so(1,n+1)\ltimes \rr^{1,n+1}$, by Corollary \ref{nabvk-cor} we have  \[\phi\in \stab(\rr v)\subset \so(1,n+1)\] for
$\phi=\nabla K$. This stabiliser is equal to the Lie algebra of similarity transformations of $\rrn$,
\be \stab(\rr v)&=&\siml(n)\\
&=&(\rr\+\so(n))\ltimes \rr^n 
=
\left\{
\left(\begin{array}{rcc}a & u^\top & 0 \\ 0& F&-u
\\0&0&-a
\end{array}\right)\ \left|\
\begin{array}{l}
a\in \rr\\
F\in\so(n)\\
u\in \rr^n
\end{array}\right.
\right\},
\ee
which is the minimal  parabolic subalgebra in $\so(1,n+1)$.
Hence we obtain

\bfolg\label{dim-cor}
The evaluation map $\kappa$ in \eqref{eval} is an injective vector space homomorphism
\belabel{evalpp}
\begin{array}{rcl}
\kappa\ : \ \k&\hookrightarrow & \siml(n)\ltimes \rr^{1,n+1}\\
K& 
\mapsto& \left(
 \left(
 \begin{array}{ccc}a& Y & 0 \\ 0&-F&-Y^\top
\\0&0&-a
\end{array}
\right), 
\left(
\begin{array}{c}
-c
\\
-X
\\
-b
\end{array}\right)
\right)
\end{array}\eeq
In particular, 
 \[1\le \dim(\k) \le 
 (2n+3)+\frac{1}{2}n(n-1)
.\]
 \efolg
 
%
% 
% Note that, when fixing coordinates of Lemma \ref{coordlemma}, by the formulae \eqref{killinit} the injection \eqref{evalsim} is explicitly given as
% \be
% \lefteqn{
% K:=-\left( ax^-+c+ \dot \Psi^\top\x \right)\del_- + \left(\Psi+\r\cdot \x\right)^i \del_i +(ax^+b) \del_+}
% \\
% &
% \longmapsto& \left(
% \left(
% \begin{array}{ccc}a & -\dot\Psi(0)^\top & 0 \\ 0& \r&\dot\Psi(0)
%\\0&0&-a
%\end{array}
%\right), 
%\left(
%\begin{array}{c}
%c
%\\
%\Psi(0)
%\\
%b
%\end{array}\right)
%\right)\ \in\  \siml(n)\ltimes \rr^{1,n+1}.
%\ee
Unfortunately, the map in \eqref{evalpp} is not a Lie algebra homomorphism. In fact, a direct computation using the bracket formula \eqref{bracket} confirms the observation \eqref{lahom}  in the general setting and yields
\belabel{lahom-exp}
\begin{array}{rcl}
[\kappa(K),\kappa(\hat K)] -\kappa([K,\hat K]) &=&\R(K,\hat K, \del_i,\del_+)|_p
\\[2mm]&=&
 \left(
 \left(
 \begin{array}{ccc}0 & (b S\hat X-\hat b S X)^\top  & 0 \\ 0& 0&\hat b S X-b S\hat X
\\0&0&0
\end{array}
\right), 
0\right)
\end{array}
,\eeq
where $S=
\hess(H)|_p$.
As a remedy, we consider the vector space
\[
\k_p(V)=\{K\in \k\mid \gg(K,V)|_p=0\}.
\]
According to Theorem \ref{killtheo}, when using coordinates of Lemma \ref{coordlemma} around $p$ elements in $\k_p(V)$ are characterised by the condition $b=0$. Hence, consulting  formula \eqref{bracket} for the Lie bracket of two Killing vector fields, we make the following observation
\bfolg
$\k_p(V)$ is a Lie subalgebra of $\k$. Moreover, the evaluation map at $p$, when restricted to $\k_p(V)$ is an injective Lie algebra homomorphism
\be
\kappa:\k_p(V)&\hookrightarrow&\co(n)\ltimes \he(n)
\\
K &
 \mapsto&
 \left(
 \begin{array}{ccc}a & -Y & c \\ 0& F&X
\\0&0&0
\end{array}
\right),
\ee
where $\co(n):=\rr\+\so(n)$ denotes the conformal Lie algebra and $\he(n)$ the $(2n+1)$- dimensional Heisenberg algebra.
\efolg
\bprf
That the evaluation map $\kappa$ at $p$ becomes a Lie algebra monomorphism follows from observation \eqref{lahom} and the defining property of pp-waves, which ensures that 
$\R(K,\hat K,\cdot,\cdot )|_p=0$ whenever $ K,\hat K\in \k_p(V)$. It can also seen 
immediately from  Theorem~\ref{killtheo}, $b=0$  or from the observation \eqref{lahom-exp}. Moreover, if $b=0$, the image of $K_p$ lies in $v^\perp$, i.e., $\kappa(\k_p(V))\subset \siml(n)\ltimes v^\perp$. Hence it remains to establish that
\be
\siml(n)\ltimes v^\perp& \simeq &\co(n)\ltimes \he(n)
\\
 \left(
 \left(
 \begin{array}{ccc}a & Y^\top & 0 \\ 0& F&-Y
\\0&0&-a
\end{array}
\right), 
\left(
\begin{array}{c}
c
\\
X
\\
0
\end{array}\right)\right)
&\mapsto
&
 \left(
 \begin{array}{ccc}a & Y^\top & c \\ 0& F&X
\\0&0&0
\end{array}
\right)
\ee
is indeed a Lie algebra isomorphism. But this is a straightforward computation.
\eprf
Because of Lemma \ref{kill-lem}, for the subalgebra
\[
\k(V)=\{K\in \k\mid \gg(K,V)=0\}
\]
of $\k_p(V)$, whose elements are characterised by $a=b=0$, we obtain
\bfolg\label{kvfolg}
The evaluation map at $p$, when restricted to $\k(V)$ is an injective Lie algebra homomorphism
\be
\kappa:\k(V)&\hookrightarrow&\so(n)\ltimes \he(n)
\\
K &
 \longmapsto&
 \left(
 \begin{array}{ccc}0& -Y & c \\ 0& F&X
\\0&0&0
\end{array}
\right).
\ee
\efolg

Returning to the evaluation map, we note that 
%$\he(n)$ contains an Abelian ideal that is also invariant under $\co(n)$. Hence
the Lie algebra $\co(n)\ltimes \he (n)$ contains an Abelian ideal
\[\a:=
\left\{	\left(\begin{array}{ccc}a & Y & c \\ 0& 0&0 \\0&0&0
\end{array}\right)\ \Biggl|\ Y\in \rrn,\ a\in\rr,\ c\in \rr\right\}\subset \co(n)\ltimes \he(n).
\]
Therefore, the quotient $(\co(n)\ltimes \he(n))/\a$ is a Lie algebra which turns out to be isomorphic to the Lie algebra of Euclidean motions $\so(n)\ltimes \rrn$ via
\be
( \co(n)\ltimes \he(n)) /\a&\simeq & \so(n)\ltimes \rrn
\\
 \begin{pmatrix}
 a & Y^\top & c \\ 0& F&X
\\0&0&0
\end{pmatrix}
+\a
&\mapsto
&
 \begin{pmatrix}
  F&X
\\0&0
\end{pmatrix}.
 \ee
Hence, we obtain
\bfolg
\label{lambda-cor}
The evaluation map $\kappa$ induces  a Lie algebra homomorphism $\lambda: \k_p(V)\to \so(n)\ltimes \rrn$ given by
\be
\k_p(V)&\stackrel{\lambda}{\longrightarrow}&\so(n)\ltimes \rrn
\\
K
 &
 \longmapsto&
 \begin{pmatrix} -F&-X
\\0&0
\end{pmatrix}.
\ee
Moreover, if $\k_p(V)$ at $p$ spans $V^\perp|_p$, then $\g:=\lambda(\k_p(V))\subset \so(n)\ltimes\rrn$ is a subalgebra that acts indecomposably on $\rr^{1,n+1}$ via
\[
\begin{pmatrix}
0&X^\top&0
\\
 0&-F&-X
\\0&0&0
\end{pmatrix}.
\]
\efolg
\bprf
Since there are Killing vector fields that span $V^\perp|_p$, by the definition of $\lambda$ 
for the projection  $\mathrm{pr}_{\rrn}:\so(n)\ltimes \rrn\to \rrn$ 
onto the translations
 we have that 
\[\mathrm{pr}_{\rrn}(\lambda(\k_p(V)))=\rrn.\]
This implies that $\g=\lambda(\k_p(V))$ acts indecomposably on $\rr^{1,n+1}$. 
\eprf
For the Killing vector fields $\k(V)$ that are tangent to $V^\perp$ we consider the ideal 
\[
\mathfrak{b}:=\rr^n\ltimes \rr \subset \so(n)\ltimes \he(n),
\]
for which we have $(\so(n)\ltimes \he(n))/\mathfrak{b}\simeq\so(n)\ltimes \rrn$.
In Corollary \ref{kvfolg}, this ideal corresponds to the elements
\[
\begin{pmatrix}
0& -Y & c \\
0& 0&0\\
0&0&0
\end{pmatrix}.
\]
In complete analogy to Corollary \ref{lambda-cor} we obtain from Corollary \ref{kvfolg} the following
\bfolg
\label{Vlambda-cor}
The evaluation map $\kappa$ induces  a Lie algebra homomorphism
$\lambda: \k(V)\to \so(n)\ltimes \rrn$.
% given by
%\be
%\k(V)&\stackrel{\lambda}{\longrightarrow}&\so(n)\ltimes \rrn
%\\
%K
% &
% \longmapsto&
% \left(
% \begin{array}{cc} -\r&-X
%\\0&0
%\end{array}
%\right).
%\ee
Moreover, if $\k(V)$  spans $V^\perp$, then $\h:=\lambda(\k(V))\subset \so(n)\ltimes\rrn$ is a subalgebra that acts indecomposably on $\rr^{1,n+1}$ as in Corollary \ref{lambda-cor}.
\efolg

For the proof of Theorem \ref{maintheo} we will  need a description of subalgebras of $\siml(n)$ that act indecomposably on $\rr^{1,n+1}$.
% Lie algebras of subgroups of Euclidean similarity transformations that act transitively on $\rrn$.
  Fortunately, there is such a  classification   due to B\'{e}rard-Bergery and Ikemakhen
  %Alekseevsky \cite{alekseevsky75, alekseevskij-vinberg-solodovnikov93}, see also
   \cite{bb-ike93}: 
%   and   \cite[Theorem 3.4]{galaev-leistner-esi}. 
\bs
%\cite{alekseevsky75, alekseevskij-vinberg-solodovnikov93}\label{gsimprop}
Let $\g\subset \siml(n)=(\rr\+\so(n))\ltimes\rrn$ act indecomposably  on $\rr^{1,n+1}$. Then either $\g$ contains the translations $\rrn$, 
%in which case $\g$ is isomorphic to  $\h\ltimes \rrn$, where $\h\subset\so(n)$ is any subalgebra, 
or  $\g$ contains $\rr^q$ for $1<q<n$, in which case there is a subalgebra $\h\subset \so(q)$  and a surjective linear map $\varphi:\h\to \rr^{n-q}$ such that  $\g$ is of the form
\belabel{gsim}
\g = \left\{ \left.
\left(
\begin{array}{cccc}
0&X&\varphi(F)&0\\
0&  F   &0&-X\\
0&0&0 &-\varphi(F)\\
0&0&0&0
\end{array}
\right)\ \right|\ F\in \h ,X\in \mathds{R}^q \right\}.
\eeq
\label{simprop}
\es
The important property in this proposition is that the rotational part $F$ of a transitively acting group of similarity transformations acts only on $\rr^p$ and annihilates the corresponding translational part $\vf(F)$.

With this at hand we are ready to prove Theorem \ref{maintheo}.
 \bprf[{\bf Proof of Theorem \ref{maintheo}}]
By the defining property \eqref{planewave} of a plane wave,
we have to show that at each point $p\in \M$ we have
$\nabla_U \R|_p=0$ for all $U\in V^\perp|_p$. 
%Since  $(\M,\gg)$ is assumed to be strongly indecomposable, by Lemma \ref{indeclemma}, without loss of generality,  we can choose $p$ such that the matrix $\R(E_+,E_i,E_+,E_j)$ is invertible.
Working with a basis of the form \eqref{basis}, from the formulae for the curvature and the Levi-Civita connection of  a pp-wave   it follows that the only possibly non-vanishing terms of $\nabla\R$  are $\nabla_{E_+}\R(E_+,E_i,E_+,E_j)$ and 
\belabel{symR}
\begin{array}{rcl}
\nabla_k\R_{ij}&:=&
\nabla_{E_k}\R(E_+,E_i,E_+,E_j)\\
&=&
\nabla_{E_i}\R(E_+,E_k,E_+,E_j)
=
\nabla_{E_j}\R(E_+,E_k,E_+,E_i),
\end{array}
\eeq
for $i,j,k=1, \ldots, n$ and,  because of the Bianchi identity, being symmetric in those.
We will now use the integrability condition \eqref{kill-int}
to show that this term also vanishes. Because of our assumption that the curvature has rank greater than $1$ almost everywhere, it suffices to work at a $p\in \M$ at which the rank of $\R$ is greater than $1$. This just means that the rank of the matrix
\[
\R_{ij}:=\R(E_+,E_i,E_+,E_j)
\]
is greater than $1$.

Since there are Killing vector fields that span $V^\perp|_p$, we can apply  Corollary \ref{lambda-cor}
and Proposition \ref{simprop} to $\g=\lambda(\k_p(V))$ giving two possible cases for $\g$.  In the first case, $\g$ contains the translations $\rrn$, i.e.,  there are Killing vector fields $K_1, \ldots, K_n$ with
\[\lambda(K_k)
=
\left(\begin{array}{ccc}
0&\e^\top_k&0
\\
0&0&-\e_k
\\
0&0&0
\end{array}
\right)\in \siml(n).\]
By the definition of $\lambda$ 
%we have that 
%\[K_i= 
%-(\lambda_ix^-+ \Psi_i^\top(x^+) )\del_- 
%
%
and recalling \eqref{atp}, this implies
 $K_k|_p=E_k$ and for the $\phi_k=\nabla K_k|_p$ that
\belabel{phirn}
\begin{array}{rcl}
\phi_k(E_j)&\in& \rr V_p, \ \text{ for }j=1, \ldots, n\\
\phi_k(E_+)&=& a_k \del_+\mod V^\perp|_p,
\end{array}\eeq
for $k=1, \ldots , n$. Without loss of generality we may assume that all but one $a_i$ are equal to zero. Indeed, the linear map from $\mathrm{span}(K_i)_{i=1}^n$ to $\rr$ defined by assigning $a_i$ to each $K_i$ has a kernel of dimension at least $n-1$. Hence, we can chose at least $n-1$ linearly independent Killing vector fields in its kernel and possible one  that is transversal to the kernel. The latter can be chosen in a way that, at $p$, it is orthogonal to the kernel, whereas the ones in the kernel can be chosen to be orthonormal to each other at $p$.

Hence, we can  assume that $a_1= \ldots =a_{n-1}=0 $, and the integrability condition \eqref{kill-int} becomes
\belabel{intrn}
\begin{array}{rcl}
\nabla_{k}\R_{ij}&=&\R(\phi_k(E_+),E_i,E_+,E_j)
+ \R(E_+,\phi_k(E_i),E_+,E_j)
\\&&{}+\R(\phi_k(E_+),E_j,E_+,E_i)+ \R(E_+,\phi_k(E_j),E_+,E_i)
\\
&=&2 a_k \R_{ij},
\end{array}
\eeq
for $i,j,k=1, \ldots, n$. Therefore, we get
\[
\nabla_{k}\R_{ij}=0,
\]
for $k=1, \ldots, n-1$ and $i,j=1, \ldots, n$, as well as
\[
2a_n\R_{ki}=\nabla_n\R_{ki}=0
\]
for all $i=1, \ldots, n$ and $k=1, \ldots , n-1$.
Hence, if $a_n$ was not zero, $\R_{nn}$ would be the only non-vanishing component of $\R_{ij}$ which contradicts the assumption that its rank is greater than one. Hence, also $a_n=0$ and therefore $\nabla_k\R_{ij}=0$ for all $i,j,k$.

This gives us an idea how to proceed in the remaining case, in which 
 $\g$ does not contain $\rrn$, but only an $\rr^N$, for $1<N<n$. Here, according to Proposition \ref{simprop}, $\g $ is of the form 
\eqref{gsim}. In the following, we will use indices $A,B,C\ldots \in \{1, \ldots , N\}$ and $b,c,d ,\ldots \in \{N+1, \ldots , n\}$ and $i,j,k\in \{1, \ldots , n\}$.
For such $\g$'s we have $N$  Killing vector fields 
such that
\[\lambda(K_A)
=
\left(\begin{array}{cccc}
0&\e^\top_A&0&0
\\
0&0&0&-\e_A
\\
0&0&0&0
\\
0&0&0&0
\end{array}
\right)\in \so(n)\ltimes \rrn,\]
with
\belabel{hirN}
\begin{array}{rcl}
K_{A}|_p&=&E_{A}
\\
\phi_{A}(E_-)&=&a_{A}\del_-\\
\phi_A(E_i)&\in& \rr V_p, \ \text{ for }j=1, \ldots, n\\
\phi_A(E_+)&=& a_A \del_+\mod V^\perp|_p,
\end{array}\eeq
and $n-N$ Killing vector fields  $K_{b}$, with
\[\lambda(K_{b})
=
\left(\begin{array}{cccc}
0&0&\e^\top_{b}&0
\\
0&\stackrel{(b)}{F}&0&0
\\
0&0&0&-\e_{b}
\\
0&0&0&0
\end{array}
\right)\in \so(n)\ltimes\rrn.
\]
Note that by Proposition \ref{simprop} all the $\stackrel{(b)}{F}\in \so(N)$ are non-zero.
By the definition of $\lambda$ and looking at \eqref{atp} this implies for $\phi_{b}=\nabla K_{b}|_p$ that
\belabel{atprot}
\begin{array}{rcl}
K_{b}|_p&=&E_{b}
\\
\phi_{b}(E_-)&=&-a_{b}\del_-\\
\phi_{b}(E_A)&=&\stackrel{(b)}{f}\hspace{-4pt}_{A}^{~B}\mod \rr V_p
\\
\phi_{b}(E_{c})&\in &\rr V_p
\\
\phi_{b}(E_+)&\in& a_{b}\del_+ \mod  V^\perp|_p.
\end{array}
\eeq
As before, without loss of generality, we can assume that $a_N$ and $a_n$ are the only $a_i$'s that are possibly non-zero.
Then we have 
\begin{eqnarray}
\label{nabAij}
\nabla_A\R_{ij} &=&2a_A\R_{ij}
\\
\label{nabbcd}
\nabla_b\R_{cd}&=&2a_b\R_{cd}
\\
\label{nabbcA}
\nabla_b\R_{cA}&=& 2a_b\R_{cA}+\stackrel{(b)}{f}\hspace{-4pt}_{A}^{~B}\R_{cB}
\\
\label{nabbAB}
\nabla_b\R_{AB}&=&2 a_b\R_{AB}+2
\stackrel{(b)}{f}\hspace{-4pt}_{(A}^{\ \ C}\R_{{B)}C}
\end{eqnarray}
With our assumption $a_1=\ldots a_{N-1}=a_{N+1}=\ldots = a_{n-1}=0$ equation \eqref{nabAij} gives
\belabel{nabAij0}
\nabla_A\R_{ij}=0, \text{ for all }A\neq N
\eeq
and thus
\belabel{nabNAj0}
a_N\R_{Aj}=0,\ \text{ for all }A\neq N.
\eeq
Similarly, equation \eqref{nabbcd} yields
\belabel{nabbcd0}
\nabla_b\R_{cd}=0, \text{ for all }b\neq n
\eeq
and thus
\belabel{nabnbc0}
a_n\R_{bc}=0,\ \text{ for all }(b,c)\neq (n,n).
\eeq
Furthermore, using the  total symmetry of $\nabla_i\R_{jk}$ we  observe that equation \eqref{nabbcA} gives
\belabel{nabAbc}
2a_A\R_{bc}= 
2a_b\R_{cA}+\stackrel{(b)}{f}\hspace{-4pt}_{A}^{~B}\R_{cB}
\eeq
and \eqref{nabbAB} yields
\belabel{nabABc}
2a_A\R_{Bc}= 
2a_c\R_{AB}+\stackrel{(c)}{f}\hspace{-4pt}_{(A}^{~D}\R_{{B)}D}.
\eeq

With all these relations, the total symmetry of $\nabla_i\R_{jk}$
implies that the only possibly non-vanishing terms of $\nabla_i\R_{jk}$ are 
\belabel{nnn}
\begin{array}{rcl}
\nabla_N\R_{NN}&=&2a_N\R_{NN}
\\
\nabla_n\R_{nn}&=&2a_n\R_{nn}
\\
\nabla_N\R_{nN}&=&a_N\R_{nN}\ =\ 
2 a_n\R_{NN}+\stackrel{(n)}{f}\hspace{-4pt}_{N}^{~C}\R_{NC}
\\
\nabla_n\R_{nN}&=& a_N\R_{nn}\ =\ 
2a_n\R_{nN}+\stackrel{(n)}{f}\hspace{-4pt}_{N}^{~B}\R_{nB}.
\end{array}
\eeq
Now we consider two cases: First assume that $a_N\neq0$. In this case equation \eqref{nabNAj0} implies that
\belabel{RAj0}\R_{Aj}=0\ \text{ for all }A\neq N
\eeq
Evaluating \eqref{nabAbc} for $A=N$ 
yields
\belabel{anN}
2a_N\R_{bc}= 
2a_b\R_{cN}+\stackrel{(b)}{f}\hspace{-4pt}_N^{~B}\R_{cB}=
2a_b\R_{cN}
=
2a_c\R_{bN}
\eeq
since $\stackrel{(b)}{F}$ is skew and hence $\stackrel{(b)}{f}\hspace{-4pt}_N^{\ N}=0$. Evaluating this for $b\neq n$ we get that 
\belabel{Rbc0}
\R_{bc}=0,\ \text{ for all $(b,c)\neq(n,n)$.}\eeq
Moreover, equation \eqref{nabABc} for $A=B=N$ for $c\neq n$  gives
\[
2a_N\R_{Nc}= \stackrel{(c)}{f}\hspace{-4pt}_{N}^{~D}\R_{ND}=0
\]
again because of \eqref{RAj0} and the skew-symmetry of $\stackrel{(c)}{F}$. So we get
\belabel{RNb0}
\R_{Nb}=0\quad \text{ for }b\neq n.
\eeq
Putting \eqref{RAj0}, \eqref{Rbc0} and \eqref{RNb0} together we get that
$\R_{NN}$, $\R_{nn}$ and $\R_{Nn}$ are the only non vanishing components of $\R_{ij}$. According to the last two equations of \eqref{nnn} they are related by
\be
a_N\R_{nN}
&=&
 a_n\R_{NN}
 \\
 a_N\R_{nn}& =&
 a_n\R_{nN}
 \ee
 This implies that $a_n\neq0$ because otherwise $\R_{NN}$ would be the only non-vanishing  component of $\R_{ij}$ which contradicts to the rank of $\R_{ij}$ being greater than one. But this implies
 \[
 a_na_N\det \begin{pmatrix} \R_{NN}&\R_{Nn}\\\R_{nN}&\R_{nn}\end{pmatrix}
 =
 0,
 \]
 which finally leads a contradiction to the rank of $\R_{ij}$ being greater than one.
 
 It remains to derive a contradiction in the case when $a_N=0$. If also $a_n=0$ we are done, so we assume  $a_n\neq0$. In this case \eqref{nabnbc0} implies that 
 \belabel{Rbc00}
\R_{bc}=0,\ \text{ for all $(b,c)\neq(n,n)$.}\eeq
Moreover \eqref{nabAbc} for $b=n$ implies that each $(\R_{cB})_{B=1}^N$ is an eigenvector of $\stackrel{(n)}{F}$.
Since $a_n\neq 0$ is real and $\stackrel{(n)}{F}$ skew, this implies that $\R_{cB}=0$ for all $c$ and $B$. 

Moreover equation  \eqref{nabABc} for $c=n$ 
becomes
\[
-2a_n\R_{AB}=\stackrel{(n)}{f}\hspace{-4pt}_{A}^{~D}\R_{BD}+ \stackrel{(n)}{f}\hspace{-4pt}_{B}^{~D}\R_{AD}
=
\stackrel{(n)}{f}\hspace{-4pt}_{A}^{~D}\R_{BD}- \stackrel{(n)}{f}\hspace{-4pt}_{D}^{~B}\R_{AD}
\]
which just 
means that the matrix $(\R_{AB})$ is an eigenvector with eigenvalue $-2a_n$ for the adjoint action of
$\stackrel{(n)}{F}\in \so(n)$ on the symmetric matrices, i.e.,
\belabel{this}
-2a_n\mathbf{\R}
=
[ \stackrel{(n)}{F},\mathbf{\R}].
\eeq
Since $ \stackrel{(n)}{F}$, when acting on symmetric matrices via the commutator, is skew-symmetric  with respect to the trace form, which, on the other hand, is positive definite on symmetric matrices, \eqref{this} implies $\R_{AB}=0$. Hence, again $\R_{nn}$ is the only non-vanishing component of $\R_{ij}$ which contradicts our assumption that the rank of the curvature endomorphism is larger than one.
This concludes the proof of Theorem~\ref{maintheo}.
\eprf
This proof and  Corollary~\ref{kvfolg} immediately give us a 
 {\em \bf proof of Theorem~\ref{kVtheo}}  when taking into account that Killing vector fields from $\k(V)$ have  $a_i=0$ for $i=1, \ldots, n$.

\bbem\label{con-rem}
Note that our proof shows that for indecomposable homogeneous pp-waves with  $1$-dimensional Lie algebra $\k/\k(V)$, we could drop the assumption on the rank of the curvature in Corollary \ref{folg2}. Indeed, if $(\M,\gg)$ is homogeneous, at each point $p$ we have, in addition to the Killing vector fields $V,K_1, \ldots , K_n$ spanning $V_p^\perp$, a  Killing vector field $\hat K$ transversal to $V_p^\perp$. In normal Brinkmann coordinates this vector field would have $b=1$ and hence, by the assumption $\dim(\k/\k(V))=1$, all the $K_i$'s would have $a_i=0$. The proof of Theorem \ref{maintheo}  then shows that $(\M,\gg)$ is a plane wave.
\ebem

%%%%%%%%%%%%%%

%%%%%%%%%%%%%%%
\begin{appendix}
\section{Normal Brinkmann coordinates for pp-waves}
Here we prove Lemma \ref{coordlemma}.
It is well known that, since a pp-wave has a parallel null vector field,  it admits local Walker coordinates \cite{walker50I}. Evaluating the curvature condition \eqref{screen-flat} in these coordinates  yields the desired form \eqref{ppcoord}. 
We will give some more detail on this, as it gives us the opportunity  to describe the coordinate freedom: By the existence of Walker coordinates, there is a $x^+$-dependent family of one-forms $ \mu=\mu_i(x^+)\d x^i$ and a $x^+$-dependent family of Riemannian metrics $\hh=h_{ij}(x^+)\d x^i\d x^j$ and a smooth function $H=H(x^+,\x) $  such that 
\belabel{brinkmann}
\gg=2\d x^+(\d x^- + H \d x^+ +\mu) + h_{ij}\d x^i\d x^j,
\end{equation}
or, more conveniently
\belabel{brinkmann1}
\gg=2\d x^+(\d x^- + H \d x^+ +\mu^\top \d\x ) + \d\x^\top \hh \d\x,
\end{equation}
where we set 
 $\x:=(x^1, \ldots , x^n)$ and slightly abuse the notation when denoting the vector $\mu=(\mu_1, \ldots , \mu_n) $ and the matrix  $\hh=(h_{ij})$ by the same symbols as the one form and the metric.
 Note that the most general coordinate transformation preserving this form is given by
\belabel{trafo}
\wt x^-=\frac{1}{a} x^-+ F(x^+,\x)
,\ \ 
\wt \x =\wt \x(x^+, \x),\ \ 
\wt x^+=a x^++b
\end{equation}
for constants $a\neq0$ and $b$, and a function $ F$ of $x^+$ and the $x^i$'s. Then, for the new ingredients $\wt H$, $\wt \mu$ and $\wt \hh$ of the metric in form \eqref{brinkmann1} 
\[
\gg=2\d\wt x^+(\d\wt x^- + \wt H \d\wt x^+ +\wt \mu^\top \d\wt \x ) + \d\wt \x^\top \wt \hh \wt \d\x
\]
we get
the relations
\belabel
{Htrafo}
\begin{array}{rcl}
H&=&a (\wt H+ \dot  F+ \wt \mu^\top \dot {\wt\x})+\tfrac{1}{2} \dot {\wt\x}^\top \wt \hh \dot {\wt\x}
\\
\mu&=& a\grad^\hh(F) + (a\wt \mu^\top +  \dot {\wt\x}^\top\wt \hh ) D(\wt\x)\\
\hh&=&  D(\wt\x)^\top \wt h D(\wt\x)
,\end{array}\end{equation}
where $\grad^\hh( F)$ denotes the gradient of $ F$ with respect to $\hh$, and $D(\wt\x)$ the Jacobian of $\wt x$ in the $x^i$ directions.

Now we turn to pp-waves.
For the curvature of  a metric in \eqref{brinkmann}
 we compute 
\belabel{curv}
\R(X,Y)Z=\R^\hh(X,Y)Z +\big( (\d^{\nabla^\hh}\dot \hh(X,Y,Z) -\tfrac{1}{2}(\nabla^\hh_Z\d\mu)(X,Y)-(\R^\hh_{X,Y}\mu)(Z)\Big)\del_-,
\end{equation}
for $X,Y,Z$ in the span of the $\del_i$'s. Paring this with $\del_i$, condition \eqref{screen-flat} shows that $\hh$ is a family of flat Riemannian metrics, and hence, by applying a transformation as in \eqref{trafo} with $ F\equiv 0$, $a=1$ and $b=0$ preserving the form of \eqref{brinkmann} but such that $h_{ij}\equiv \delta_{ij}$.
In these coordinates, pairing \eqref{curv} with $\del_+$, condition \eqref{screen-flat} becomes 
\[
0=\nabla^\hh_{\del_i}\mu(\del_j,\del_k) = \del_i(m_{jk}),
\]
where $\d\mu = M_{ij}\d x^i \wedge \d x^j$, where the $\d$ denotes the differential only in the $x^i$-directions. Hence $M(x^+):=(M_{ij}(x^+))_{i,j=1}^n\in \so(n)$ is an $x^+$-dependent family of skew-symmetric matrices. 
For this $M$ , we consider the linear  ODE
\begin{equation}\label{ODE}
\dot A=-AM,\quad
A(0)=A_-\in \O(n).
\end{equation}
This has a unique solution $A(x^+)$ which satisfies
\[
\frac{\d}{\d x^+}(A A^\top )= -A M A^\top- A M^\top A^\top=0,
\]
 since $M$ skew. Hence, $A(0)\in \O(n)$ implies that $A(x^+)\in \O(n)$ for all $x^+$. For such a solution $A$, we define the $x^+$-dependent  one-form
 \[
 \alpha = \x^\top \dot A^\top A \d\x
 = x^l\dot A^i_l\delta_{ij}A^j_k\d x^k.
 \]
 This form $\alpha+\mu$  is closed,
\be
\d(\mu-\alpha) &= &
 M_{lk}\d x^l\wedge \d x^k 
 -
 \dot A^i_l\delta_{ij}A^j_kx^l\wedge \d x^k
\\
& =&
(M_{lk}- (\dot A^\top A)_{lk}) \d x^l\wedge \d x^k
\\
&=&
(  M_{lk}+ (  M^\top A^\top A)_{lk}) \d x^l\wedge \d x^k
\\
&=&0.
\ee
Now, for given $\mu$ in \eqref{brinkmann} wit $\hh_{ij}\equiv \delta_{ij}$, let $ F= F(x^+,\x) $ be a solution to
$\d F =\mu-\alpha $ and $A$ a solution to \eqref{ODE} and
consider the coordinate transformation
\belabel{trafo1}
\wt x^-= x^-+ F(x^+,\x)
,\ \ 
\wt\x =A\x, \text{ i.e., } \wt x^i=A^i_kx^k,\ \ 
\wt x^+=x^+.
\end{equation}
Then, according to \eqref{Htrafo}, we have $\wt \hh=\delta_{ij}$ and moreover,  
\[
\mu=
\d F + (\wt \mu^\top +  \x^\top \dot A^\top ) A\d\x
=
\d F + \wt \mu^\top   A\d\x+\alpha.
\]
Since $\d F=\mu-\alpha$, this  implies $\wt \mu=0$, as required.
Note that the general  transformation preserving the form \eqref{ppcoord} of Brinkmann coordinates are of the form
\[
\wt x^-= \frac{1}{a}x^-+ F(x^+,\x)
,\quad
\wt\x =A\x+\vect{c}(x^+),\quad 
\wt x^+=ax^++b,
\]
where $a\neq0$ and $b$ are constants, $\vect{c}(x^+)\in \rrn$, $A=A(x^+)\in \O(n)$ satisfying the PDE
\belabel{trafopde}
0= a\d F +(\dot A\x+\dot{\vect{c}})^\top A \d\x.
\end{equation}
The integrability condition for this is 
\[
0= \d\x^\top \dot A^\top A \d\x,
\]
which implies that $\dot A=0$
(note that $\d\x^\top \dot A^\top A \d\x$ is indeed a two-form,
as $\dot A^\top A$ is skew-symmetric). This implies that $ F$ is linear in the $x^i$'s, i.e.,
\[
 F(x^+,\x) = -\frac{1}{a}\dot{\vect{c}}^\top (x^+) A \x +\beta (x^+)
\] for a function $\beta=\beta(x^+)$.
Hence, 
the general  transformation preserving the form \eqref{ppcoord} of a Brinkmann coordinates are given by a constant matrix $A\in \O(n)$,  a vectorial function   $\vect{c}$ of $x^+$ and a real function $\beta$ of $x^+$, and two real numbers $a\neq 0$ and $b$, and the transformation is
\belabel{trafo1pp}
\wt x^-= \frac{1}{a}( x^- -\dot{\vect{c}}(x^+)^\top A \x )+\beta (x^+)
,\quad
\wt\x =A\x+\vect{c}(x^+),\quad
\wt x^+=ax^++b,
\end{equation} 
The function $\wt H$ is then given as
\belabel{Htrafopp}
\wt H= \frac{1}{a}(H +
\ddot{\vect{c}}(x^+)^\top  A \x)  +\dot \beta 
-
\tfrac{1}{2a} \dot{\vect{c}}(x^+)^\top \dot{\vect{c}}(x^+)
\end{equation}

Clearly, by applying a translation we can choose these 
coordinates in a way that $p$ goes to the origin.

It remains to show that for a given Brinkmann coordinates $\vf=(x^+,x^-,\x)$ mapping $p$ to the origin, there is a coordinate transformation of the form \eqref{trafo1pp}
 that fixes the origin and provides us with normal Brinkmann coordinates, i.e., for  which the new function $\tilde{H}$ satisfies
 \belabel{nbc}
 \begin{array}{rcl}
 \tilde{H}|_{\tilde{\vf}^{-1}(x^+,\oo)}&=&0
\\
 \frac{\del \tilde H}{\del\tilde{x}^i}|_{\tilde{\vf}^{-1}(x^+,\oo)}&=&0
 \end{array}
 \eeq
 for all $x^+$.
 To this end we consider a transformation \eqref{trafo1pp} with $A=\delta_{ij}$, $b=0$ and $a=1$. 
 Let $\vect{c}=(c_1, \ldots c_n)$ the solution to the ODE system
 \[
 \ddot c_i(t) =-\frac{\del}{\del x^i}H(\vf^{-1}( t,-\vect{c} (t)),
 \]
 for $i=1, \ldots, n$ with one initial condition $c_i(0)=0$. Given such a solution $\vect{c}=(c_1,\ldots , c_n)$, let $\beta$ be the solution to the ODE
 \[
\dot\beta= \frac{1}{2}\dot{\vect{c}}^\top\dot{\vect{c}} -H(\vf^{-1}( t,-\vect{c} (t)),
\]
 with the initial condition $\beta(0)=0$.
Using these solutions $\vect{c}$ and $\beta$ in the coordinate transformation \eqref{trafo1pp}, the formula \eqref{Htrafopp} shows that in the new coordinates we have  
equations \eqref{nbc}
for all $x^+$.
\end{appendix}

%\bibliographystyle{abbrv}
%%\bibliographystyle{amsalpha}
%\bibliography{GEOBIB}

\begin{thebibliography}{10}

\bibitem{AichelburgBalasin96}
P.~C. Aichelburg and H.~Balasin.
\newblock Symmetries of pp-waves with distributional profile.
\newblock {\em Classical Quantum Gravity}, 13(4):723--729, 1996.

\bibitem{AichelburgBalasin97}
P.~C. Aichelburg and H.~Balasin.
\newblock Generalized symmetries of impulsive gravitational waves.
\newblock {\em Classical Quantum Gravity}, 14(1A):A31--A41, 1997.
\newblock Geometry and physics.

\bibitem{alekseevsky-kimelfeld75}
D.~V. Alekseevski{\u\i} and B.~N. Kimel{\cprime}fel{\cprime}d.
\newblock Structure of homogeneous {R}iemannian spaces with zero {R}icci
  curvature.
\newblock {\em Funkcional. Anal. i Prilo\v zen.}, 9(2):5--11, 1975.

\bibitem{bb-ike93}
L.~B{\'e}rard-Bergery and A.~Ikemakhen.
\newblock On the holonomy of {L}orentzian manifolds.
\newblock In {\em Differential Geometry: Geometry in Mathematical Physics and
  Related Topics (Los Angeles, CA, 1990)}, volume~54 of {\em Proc. Sympos. Pure
  Math.}, pages 27--40. Amer. Math. Soc., Providence, RI, 1993.

\bibitem{blau-oloughlin03}
M.~Blau and M.~O'Loughlin.
\newblock Homogeneous plane waves.
\newblock {\em Nuclear Phys. B}, 654(1-2):135--176, 2003.

\bibitem{BondiPiraniRobinson59}
H.~Bondi, F.~A.~E. Pirani, and I.~Robinson.
\newblock Gravitational waves in general relativity. {III}. {E}xact plane
  waves.
\newblock {\em Proc. Roy. Soc. London Ser. A}, 251:519--533, 1959.

\bibitem{brinkmann25}
H.~W. Brinkmann.
\newblock Einstein spaces which are mapped conformally on each other.
\newblock {\em Math. Ann.}, 94(1):119--145, 1925.

\bibitem{cahen-wallach70}
M.~Cahen and N.~Wallach.
\newblock {L}orentzian symmetric spaces.
\newblock {\em Bull. Amer. Math. Soc.}, 79:585--591, 1970.

\bibitem{candela-flores-sanchez03}
A.~M. Candela, J.~L. Flores, and M.~S{\'a}nchez.
\newblock On general plane fronted waves. {G}eodesics.
\newblock {\em Gen. Relativity Gravitation}, 35(4):631--649, 2003.

\bibitem{ehlers-kundt62}
J.~Ehlers and W.~Kundt.
\newblock Exact solutions of the gravitational field equations.
\newblock In L.~Witten, editor, {\em Gravitation: {A}n introduction to current
  research}, pages 49--101. Wiley, New York, 1962.

\bibitem{EinsteinRosen37}
A.~Einstein and N.~Rosen.
\newblock On gravitational waves.
\newblock {\em Journal of the Franklin Institute}, 223(1):43--54, 1937.

\bibitem{Figueroa-OFarriHustler12}
J.~Figueroa-O'Farrill and N.~Hustler.
\newblock The homogeneity theorem for supergravity backgrounds.
\newblock {\em J. High Energy Phys.}, (10):014, front matter + 8, 2012.

\bibitem{Figueroa-OFarriHustler14}
J.~{Figueroa-O'Farrill} and N.~{Hustler}.
\newblock {The homogeneity theorem for supergravity backgrounds II: the
  six-dimensional theories}.
\newblock {\em Journal of High Energy Physics}, 4:131, Apr. 2014.

\bibitem{Figueroa-OFarriMeessenPhilip05}
J.~Figueroa-O'Farrill, P.~Meessen, and S.~Philip.
\newblock Supersymmetry and homogeneity of {M}-theory backgrounds.
\newblock {\em Classical Quantum Gravity}, 22(1):207--226, 2005.

\bibitem{Figueroa-OFarriPhilipMeessen05}
J.~Figueroa-O'Farrill, S.~Philip, and P.~Meessen.
\newblock Homogeneity and plane-wave limits.
\newblock {\em J. High Energy Phys.}, (5):050, 42, 2005.

\bibitem{hull84pp}
C.~M. Hull.
\newblock Exact pp-wave solutions of {$11$}-dimensional supergravity.
\newblock {\em Phys. Lett. B}, 139(1-2):39--41, 1984.

\bibitem{jek60}
P.~Jordan, J.~Ehlers, and W.~Kundt.
\newblock Strenge {L}{{\"o}}sungen der {F}eldgleichungen der allgemeinen
  {R}elativit{\"a}tstheorie.
\newblock {\em Akad. Wiss. Mainz. Abh. Math.-Nat. Kl.}, 1960:21--105, 1960.

\bibitem{jek60engl}
P.~Jordan, J.~Ehlers, and W.~Kundt.
\newblock Republication of: {E}xact solutions of the field equations of the
  general theory of relativity.
\newblock {\em Gen. Relativity Gravitation}, 41(9):2191--2280, 2009.
\newblock Translated from the German [Akad. Wiss. Mainz. Abh. Math.-Nat. Kl.
  {{\bf{1}}960}, 21--150].

\bibitem{Kostant55}
B.~Kostant.
\newblock Holonomy and the {L}ie algebra of infinitesimal motions of a
  {R}iemannian manifold.
\newblock {\em Trans. Amer. Math. Soc.}, 80:528--542, 1955.

\bibitem{leistnerjdg}
T.~Leistner.
\newblock On the classification of {L}orentzian holonomy groups.
\newblock {\em J. Differential Geom.}, 76(3):423--484, 2007.

\bibitem{leistner-schliebner13}
T.~{Leistner} and D.~{Schliebner}.
\newblock {Completeness of compact Lorentzian manifolds with special holonomy}.
\newblock {\em Preprint}, arXiv:1306.0120, June 2013.

\bibitem{penrose76}
R.~Penrose.
\newblock Any space-time has a plane wave as a limit.
\newblock In {\em Differential geometry and relativity}, pages 271--275.
  Mathematical Phys. and Appl. Math., Vol. 3. Reidel, Dordrecht, 1976.

\bibitem{Philip06}
S.~Philip.
\newblock Penrose limits of homogeneous spaces.
\newblock {\em J. Geom. Phys.}, 56(9):1516--1533, 2006.

\bibitem{sippel-goenner86}
R.~Sippel and H.~Goenner.
\newblock Symmetry classes of pp-waves.
\newblock {\em Gen. Relativity Gravitation}, 18(12):1229--1243, 1986.

\bibitem{walker50I}
A.~G. Walker.
\newblock Canonical form for a {R}iemannian space with a parallel field of null
  planes.
\newblock {\em Quart. J. Math., Oxford Ser. (2)}, 1:69--79, 1950.

\end{thebibliography}
%
%\end{document}
\def\cprime{$'$} \def\cprime{$'$} \def\cprime{$'$}

\end{document}